\documentclass[11pt,twoside]{article}
\usepackage{amsmath,amsfonts,amssymb,cases,mathdots,color,colordvi,url,tikz}
\usepackage{hyperref}
\usepackage{bm}
\usepackage{bbm}
\usepackage{pgfpages}
\usepackage{tabularx}
\usepackage{graphics}%\textwidth 16.cm
\usepackage{subfigure}
\usepackage{epsfig}
\usepackage[figuresright]{rotating}
\usepackage{rotating}
\usepackage{algorithm}
\usepackage{algorithmic}
\usepackage{physics}
\usepackage{epstopdf}
\usepackage{parskip}
\usepackage{footmisc}
\usepackage{tablefootnote}
\usepackage[justification=centering]{caption}
\usepackage[normalem]{ulem}

\newtheorem{proposition}{Proposition}%[section]
\newtheorem{theorem}[proposition]{Theorem}
\newtheorem{lemma}[proposition]{Lemma}

\newtheorem{definition}[proposition]{Definition}

\newtheorem{remark}{Remark}
\newtheorem{assumption}{Assumption}

\newcommand{\be}{\begin{equation}}
	\newcommand{\ee}{\end{equation}}
\newcommand{\ba}{\begin{eqnarray}}
	\newcommand{\ea}{\end{eqnarray}}
\newcommand{\bas}{\begin{eqnarray*}}
	\newcommand{\eas}{\end{eqnarray*}}

\graphicspath{ {./Figures/} }
\def\H{{\mathcal{H}}}

\def\mbR{{\mathbb R}}
\def\Re{{\mathbb R}}

\def\S{{\mathcal S}}
\def\C{{\mathcal C}}

\def\F{{\mathcal F}}
\def\G{{\mathcal G}}
\def\I{{\mathcal I}}
\def\J{{\mathcal J}}

\def\cL{{\mathcal L}}

\def\N{{\mathcal N}}
\def\O{{\mathcal O}}

\def\U{{\mathcal U}}
\def\V{{\mathcal V}}

\def\bfq{{\bf q}}

\def\bfs{{\bf s}}

\def\bfa{{\bf a}}
\def\bfb{{\bf b}}
\def\bfc{{\bf c}}
\def\bfd{{\bf d}}

\def\bft{{\bf t}}
\def\bfu{{\bf u}}
\def\bfv{{\bf v}}
\def\bfx{{\bf x}}
\def\bfy{{\bf y}}
\def\bfz{{\bf z}}
\def\bfw{{\bf w}}

\def\bfog{{\bf {\boldsymbol \omega}}}
\def\bfone{{\bf 1}}

\def\OI{{\overline{\I}}}
\def\OJ{{\overline{\J}}}

\def\dom{{\rm dom}\;}
\def\ri{{\rm ri}}

\def\why{\widehat{\bfy}}
\def\whz{\widehat{\bfz}}
\def\whw{\widehat{\bfw}}
\def\argmin{\mathop{{\rm argmin}}}

\def\obx{\overline{\bfx}}

\def\bp{ \textbf{Proof.} }
\def\ep{ \hfill $\Box$ }

\def\oT{ \overline{T} }
\def\omS{\overline{\S}}

\def\bfbeta{\boldsymbol{\beta}}
\def\bfepsilon{\boldsymbol{\epsilon}}
\def\bfla{{\boldsymbol{\lambda}}}

\def\bfga{\boldsymbol{\gamma}}
\def\whp{\widehat{\partial}}

\def\bfmu{{\boldsymbol{\mu}}}

\def\bfzt{{\bf {\boldsymbol \zeta}}}

\def\Rge{\mbox{\rm Rge}\,}
\def\Ker{\mbox{\rm Ker}\,}
\def\lev{\mbox{\rm lev}}
\def\conv{\mbox{\rm conv}}
\def\Diag{\mbox{Diag}}
\def\sgn{\mbox{sgn}}
\def\mbone{\mathbbm{1}}

\def\wtB{\widetilde{B}}
\def\wtb{\widetilde{\bfb}}

\title{On the Stationary Duality of Structural Composite Cardinality Optimization}
\author{Penghe Zhang\thanks{Department of Data Science and Artificial Intelligence, The Hong Kong Polytechnic University, Hong Kong SAR, China, E-mail: {penghe.zhang@polyu.edu.hk} },
\ \ Naihua Xiu\thanks{School of Mathematics and Statistics, Beijing Jiaotong University, Beijing 100044, China, E-mail: {nhxiu@bjtu.edu.cn} } \ \ and \ \
Houduo Qi\thanks{Department of Data Science and Artificial Intelligence, and Department of Applied Mathematics, The Hong Kong Polytechnic University, Hong Kong SAR, China, E-mail: {houduo.qi@polyu.edu.hk} }
	}

\date{}

\begin{document}
	\maketitle	

	\begin{abstract}
Simple cardinality refers to counting nonzero elements of an independent variable
satisfying certain properties. 
Composite cardinality is a simple counting process composited with an affine mapping, and is therefore more complicated than the simple cardinality.
We study the composite cardinality optimization problem (CCOP) with structures covering a wide range of applications. Through the use of the stationary duality, we reduce the composite counting to simple counting, and thereby obtain a dual formulation of CCOP. For both primal and dual problems, we investigate the sufficient conditions for the existence of global
solutions. Those conditions are validated on representative examples from existing literature. We then show that local solutions of the primal and dual problems are equivalent to their stationary points. This result further helps us establish a one-to-one correspondences between primal and dual local solutions. 
We also demonstrate that the correspondence holds for a pair of global solutions to the primal and dual problems, provided that the dual weighted parameters are appropriately selected. The reported theoretical results lay foundation for developing numerical algorithms for CCOP in future.

	\vspace{3mm}
	
	
	%\noindent{\bf \textbf{Mathematical Subject Classification}:} 90C26, 90C30, 90C90
\end{abstract}
{}
%\tableofcontents
%\numberwithin{equation}{section}

%%%%%%%%%%%%%%%%%%%%%%%%%%%%%%%%%%%%%%%%%%%%%
\section{Introduction}

This paper is concerned with the following Composite Cardinality Optimization
Problem (CCOP):
\begin{align} \tag{P} \label{P}
	\inf_{\bfx \in \mbR^n} \; F(\bfx) := f(\bfx) + g(A\bfx) + \Phi_\bfla(B \bfx - \bfb), 
\end{align}
where ``$:=$'' means ``define'', $f:\mbR^n \to (-\infty, \infty]$ and $g:\mbR^m \to (-\infty, \infty]$ are proper, lower semi-continuous (lsc), and convex functions, 
and the cardinality function $\Phi_\bfla (\bfu)$ counts the nonzero elements in $\bfu$ satisfying certain properties.
Moreover, $A \in \Re^{m \times n}$, $B \in \Re^{r \times n}$,
$\bfb \in \Re^r$, and $\Phi_\bfla \in \{ \Phi_{0,\bfla}, \Phi_{+,\bfla}\}$. For given vector $\bfu = (u_1, \cdots, u_r)^\top$, we define
%Typical choices for $\Phi_\bfla$ are
%\Phi_{0,\bfla}^\bfla (\bfu)
\[
\Phi_{0,\bfla}(\bfu) := \sum_{i = 1}^r \lambda_i \mbone_{\{ u_i \neq 0 \}} \quad \mbox{and} \quad
\Phi_{+,\bfla} (\bfu)  := \sum_{i = 1}^r \lambda_i \mbone_{\{ u_i > 0 \}} ,
\]
where $\bfla : = (\lambda_1, \cdots, \lambda_r)^\top >0$ is the given regularization parameter vector,
and $\mbone_{\{ \O \}}$ is the characteristic function that returns 1 when the formula $\O$ holds, otherwise returns 0. 
In particular, when each component of $\bfla$ equals 1, $\Phi_{0,\bfla}$ reduces to the $\ell_0$ norm $\| \cdot \|_0$ and $\Phi_{+,\bfla}$ becomes 
the sum of $0/1$-loss terms (aka Heaviside function). 
Problem \eqref{P} consists of three functions and is often referred to as a three-block problem. If $\inf_{\bfx \in \mbR^n} F(\bfx)$ is finite and attainable at a global solution, we replace it with $\min_{\bfx \in \mbR^n} F(\bfx)$. Later, we will provide sufficient conditions to ensure the existence of global solutions to \eqref{P}.

The purpose of this paper is to conduct a comprehensive study
on the {\em stationary dual} problem
(in the form of minimization):
\begin{align}  \tag{D} \label{D}
	\inf_{\bfw = [\bfy; \bfz] \in \mbR^m \times \mbR^r} \
	G(\bfw) := \langle \bfb, \bfz \rangle 
	+ f^*(-A^\top \bfy -B^\top \bfz) + g^*(\bfy) + \Psi_\bfmu(\bfz),
\end{align}
where $f^*$ and $g^*$ are respectively the conjugate functions of
$f$ and $g$ in convex analysis, and $\Psi_\bfmu \in \{ \Psi_{0,\bfmu}, \Psi_{+,\bfmu}\}$. For given vector $\bfz = (z_1,\cdots,z_r)^\top$, we define
\begin{align*}
	&\Psi_{0,\bfmu}(\bfz) := \sum_{i = 1}^r \mu_i \mbone_{\{ z_i \neq 0 \}} \ \ && \mbox{(corresponding to $\Phi_\bfla= \Phi_{0,\bfla}$)} \notag \\
	& \Psi_{+,\bfmu}(\bfz) := \sum_{i = 1}^r \mu_i \mbone_{\{ z_i \neq 0 \}} + \delta(\bfz|\Re^r_+) \ &&\mbox{(corresponding to $\Phi_\bfla=\Phi_{+,\bfla}$)}. \notag
\end{align*}
Here $\bfmu:=(\mu_1,\cdots,\mu_r)^\top > 0$ is the regularization parameter vector, and for a given set $\U$, $\delta(\bfx|\U)$ is the indicator function which is $0$ whenever $\bfx \in \U$ and $+\infty$ otherwise.
We postpone the derivation of the dual problem to avoid heavy 
calculation at this stage. 

The main reason for considering the dual problem is due to
computational concern. 
The composition form of the cardinality function
$\Phi_\bfla$ with the linear
mapping $B \not= I$ (not identity matrix) makes the counting process significantly more complex. For example, $\Phi_{0,\bfla}(\bfu)$ admits closed-form proximal operator, whereas $\Phi_{0,\bfla}(B\bfx - \bfb)$ does not for an arbitrary matrix $B$. This poses challenges for both algorithmic development and convergence analysis. For primal problem with a composite term $\Phi_\bfla(B\bfx - \bfb)$, the majority of existing continuous optimization methods require a regularity condition on matrix $B$ (e.g. surjectivity) to ensure global convergence (see e.g. \cite{bolte2018nonconvex,bot2019proximal,li2015global,zhou2021quadratic}).

%
%For example, the region 
%$ \Omega :=\Re^r_+ \cap \{ B\bfx - \bfb \; | \; \bfx \in \Re^r\}$
%though is polyhedral, is more complicated than $\Re^r_+$ itself.
In contrast, the dual problem reduces the composite counting $\Phi_\bfla(B\bfx -\bfb)$ to 
simple counting via $\Psi_\bfmu(\bfz)$. Therefore, the regularity condition on $B$ can be dropped in the global convergence analysis of the algorithm proposed in \cite{zhang2025composite}. However, before leveraging the computational benefits offered by model \eqref{D}, two theoretical questions must be answered: the existence of global solutions and correspondence between primal-dual solutions. Addressing the two questions constitute the main contributions of this paper.

%However, there are two basic questions to be investigate before using \eqref{D}: the existence of global solution and the correspondence between primal-dual solutions require further investigate  
%However, to consolidate the use of dual model \eqref{D}, the existence of global solution and the correspondence between primal-dual solutions require further investigate.
%\blue{however, the existence of solution, correspondence of primal-dual solutions}

In the following, we first show that the primal problem
\eqref{P} covers a wide range of applications via concrete examples.
We then discuss the major research questions with the relevant 
references. 
Finally, we summarize our main contributions on those questions.

%This has made the counting over $\Omega$ significantly difficult than
%that over $\Re^r_+$ itself.

%As we will review below, the model is quite general and covers a wide range 
%of applications. The often studied case is $B=I$ is the identity matrix and $\bfb=0$.
%However, the composition nature of the cardinality function with the linear
%mapping $B \not= I$ makes the counting significantly more complex. For example, the region 
%$ \Omega :=\Re^r_+ \cap \{ B\bfx - \bfb \; | \; \bfx \in \Re^r\}$
%though is polyhedral, is more complicated than $\Re^r_+$ itself.
%This has made the counting over $\Omega$ significantly difficult than
%that over $\Re^r_+$ itself.
%The main purpose of this paper is to extend the recently proposed stationary
%duality in \cite{zhang2025composite}	
%to \eqref{P} and the dual problem gets around this difficulty.
%In the following, we first list a few important applications modeled by 
%\eqref{P}. They will be used to illustrate our theoretical results. 
%We then state the dual problem, especially highlighting its root in 
%convex duality. Finally, we will explain our main contribution in terms
%of the research questions related to the dual problem.
%Relevant references and existing results will be reviewed throughout.

%%-----------------------------------
\subsection{Examples} \label{sec-examples}

In this part, we show that the three-block structure in \eqref{P} is 
adequate in covering a wide range of applications. 
We particularly list two of them for their use in the demonstration of our obtained theoretical results.

{\bf (E1) Support Vector Machine (SVM) with Heaviside loss} \cite{vapnik1999nature,zhang2025composite}.
Suppose we have $r$ data points $\{ \bfx_i, c_i\}$, $i\in [r]$ with each $\bfx_i \in \Re^s$ being the feature vector and
$c_i \in \{-1,1\}$ being its class label. 
The SVM is to construct an optimal hyperplane $\{ \bfx \; | \; \langle\bfx, \bfog \rangle + \omega_0 =0\}$ with $\bfog \in \Re^s$ and $\omega_0 \in \Re$ being a bias to separate the data into two classes according to their respective labels. 
In practice, the data are usually not perfectly classified.
Therefore, it is desired to construct an optimal plane that has the least number of mis-classified points. This principle leads to the following formulation:
\[
\min_{\bfog, \omega_0} \; \frac 12 \| \bfog \|^2 + \lambda \sum_{i=1}^r h \Big( 1 - c_i (\langle\bfog, \bfx_i \rangle + \omega_0 ) \Big),
\]
where $\lambda > 0$, and $h(t) = 1$ for $t>0$ and $0$ otherwise (known as Heaviside function or $0/1$-loss). 
Let 
\[
Q := \left[
\begin{array}{c}
	\bfq_1^\top \\
	\vdots  \\
	\bfq_r^\top
\end{array} 
\right], \quad 
\overline{Q} := \left[
\begin{array}{cc}
	Q, & \bfone_r 
\end{array} 
\right], 
\quad
\bfc := \left[
\begin{array}{c}
	c_1 \\
	\vdots \\
	c_r
\end{array} 
\right],
\quad
\bfx := \left[
\begin{array}{c}
	\bfog \\
	\omega_0
\end{array} 
\right].
\]
Then the SVM model can be represented as
\be \label{hSVM} \tag{Heaviside-loss SVM}
\min_{\bfx= [\bfog;\omega_0]} \ \frac 12 \| \bfog\|^2 + \Phi_{+,\bfla}\Big(\bfone_r - \Diag(\bfc) \overline{Q} \bfx \Big),
\ee 
which in the form of \eqref{P} has $f(\bfx) = \| \bfog\|^2/2$, $g \equiv 0$, $\Phi_\bfla(\bfu) = \Phi_{+,\bfla}(\bfu)$, $\bfla = \lambda \bfone_r$, $B=- \Diag(\bfc)\overline{Q}$, and $\bfb = -\bfone_r$. Here $\bfone_r$ is the $r$-dimensional vector with each component being 1.

%{\bf (E2) The $\ell_0$-edge-denoising model} \cite{fan2018approximate}. 

%Variations of this model have appeared in many applications. For example, in calcium
%imaging, the task is to determine the exact moment in time at which a neuron spikes.
%In this case, $(V, E)$ is a line-graph with $(i, i+1) \in E$ for $i \in [n-1]$ ($n$ vertices and $(n-1)$ edges). 
%The model assumes that the observed fluorescence is a noisy version of the underlying
%calcium, which exponentially decay, unless there is a spike. 
%We copy the optimization formulation \cite[Eq.~(4)]{jewell2018exact} below (notation: here we used $\bfbeta$ for $\bfc$ and $\bfb$ for $\bfy$)
%\be \label{Calcium} \tag{Calcium-Imaging}
%\min_{\bfbeta \ge 0} \; \frac 12 \| \bfbeta - \bfb\|^2 + \lambda \| (D\bfbeta)_+\|_0,
%\ee 
%where $(D\bfbeta)_{(i, i+1)} = \beta_{i+1} - \gamma \beta_i$ for $i \in [n-1]$ and
% $\gamma>0$ is a decaying factor. In \cite{jewell2018exact}, the constraint $\bfbeta \ge 0$ is implicitly implied. 
%We may cast this problem into the format of \eqref{P} with the linear operator $A$ setting to be the identity matrix:
% \[
%  f(\bfbeta) = \delta(\bfbeta \; | \; \Re^n_+),\quad
%  g(\bfbeta) = \frac 12\| \bfbeta - \bfb\|^2, \quad
%  \mbox{and} \quad
%  \Phi_\bfla(\bfbeta) = \| (D \bfbeta)_+\|_0 .
% \]
% This also belongs to the problem of {\em multiple change-point detection} listed as
% the first applications in \cite{fan2018approximate}. 

{\bf (E2) Sparsity-driven energy minimization} \cite{nikolova2011energy}. 
This is a general model and it comes under various names in different applications, for instance, image deblurring, segmentation and reconstruction in inverse problems.
It can be stated as follows.
\be \label{Energy-minimization} \tag{Energy-minimization}
\min_{\bfx \in \Omega} \; \frac 12 \| A\bfx - \bfa \|^2 + \Phi_{0,\bfla}(D\bfx),
\ee 
where $A \in \Re^n \to \Re^m$ and $D: \Re^n \to \Re^r$ be linear operators (e.g., a convolution operator and a finite-difference or wavelet-frame operator, respectively), $\bfa \in \Re^m$ is given
observation and $\bfla > 0$ is a weighted parameter vector (e.g. \cite{zhang2013}). Furthermore, $\Omega \subseteq \Re^n$ is a closed convex set such as a non-negative constraint $ \{ \bfx | \bfx \ge 0\}$ or a box constraint $\{ \bfx | L \leq \bfx \leq U \}$.
The least-square term can be cast as an energy function. The additional constraint
$\bfx \in \Omega$ is also a distinctive feature and can cover a range of applications.
Obviously, %all the three models can be put in the form of \eqref{P}. For instance, 
\eqref{Energy-minimization} can be cast in the form \eqref{P}
\[
f(\bfx) = \delta (\bfx \; | \; \Omega), \quad
g(\bfv) = \frac 12 \| \bfv - \bfa\|^2, \quad
\Phi_\bfla(\bfu) = \Phi_{0,\bfla}(\bfu),   \ B= D, \ \mbox{and} \ \bfb=0.
\]

The \eqref{Energy-minimization} model also covers the $\ell_0$-edge-denoising problem
\cite{fan2018approximate}.
We briefly describe it below.
Suppose $G=(V, E)$ be a fully connected (undirected) graph with vertices $V:=[n]=\{1, \ldots, n\}$
and the edge set $E$.
Let $\widehat{\bfbeta} = (\widehat{\beta}_1, \cdots, \widehat{\beta}_n) \in \Re^n$ be such that its $i$th value (e.g., unknown signal value) is associated with the vertex $i$. It is observed in $\bfbeta = \widehat{\bfbeta} + \bfepsilon$ with
$\bfepsilon$ be i.i.d Gaussian noises with $0$ mean. 
The purpose is to recover the true signal $\widehat{\bfbeta}$ through the observations in $\bfbeta$.
When $\widehat{\bfbeta}$ is a piecewise signal over the graph, the set of edges ${i,j} \in E$ satisfying $\widehat{\beta}_{i} \not=\widehat{\beta}_{j}$ is a small subset of $E$. This has led to the $\ell_0$-edge-denoising model:
\be \label{Edge-Denoising} \tag{Edge-Denoising}
\min_{\bfx} \; \frac 12 \| \bfx - \bfbeta \|^2 + \Phi_{0,\bfla}(D\bfx),
\ee 
where $D: \Re^n \to \Re^{|E|}$ satisfies $(D \bfx)_{(i,j)} = x_j - x_i$ for $(i,j) \in E$. 
In the particular case where $(V, E)$ is a line-graph with $(i, i+1) \in E$ for $i \in [n-1]$ ($n$ vertices and $(n-1)$ edges),  the Calcium imgaging problem  \cite[Eq.~(4)]{jewell2018exact} takes form:
\be \label{Calcium} \tag{Calcium-Imaging}
\min_{\bfx \ge 0} \; \frac 12 \| \bfx - \bfbeta\|^2 + \Phi_{+,\bfla}(D\bfx),
\ee 
We omit the detailed background behind this model.

%%%%%%
%Variations of this model have appeared in many applications. For example, in calcium
%imaging, the task is to determine the exact moment in time at which a neuron spikes.
%In this case, $(V, E)$ is a line-graph with $(i, i+1) \in E$ for $i \in [n-1]$ ($n$ vertices and $(n-1)$ edges). 
%The model assumes that the observed fluorescence is a noisy version of the underlying
%calcium, which exponentially decay, unless there is a spike. 
%We copy the optimization formulation \cite[Eq.~(4)]{jewell2018exact} below (notation: here we used $\bfbeta$ for $\bfc$ and $\bfb$ for $\bfy$)
%\be \label{Calcium} \tag{Calcium-Imaging}
%\min_{\bfbeta \ge 0} \; \frac 12 \| \bfbeta - \bfb\|^2 + \lambda \| (D\bfbeta)_+\|_0,
%\ee 
%where $(D\bfbeta)_{(i, i+1)} = \beta_{i+1} - \gamma \beta_i$ for $i \in [n-1]$ and
% $\gamma>0$ is a decaying factor. In \cite{jewell2018exact}, the constraint $\bfbeta \ge 0$ is implicitly implied. 
%We may cast this problem into the format of \eqref{P} with the linear operator $A$ setting to be the identity matrix:
% \[
%  f(\bfbeta) = \delta(\bfbeta \; | \; \Re^n_+),\quad
%  g(\bfbeta) = \frac 12\| \bfbeta - \bfb\|^2, \quad
%  \mbox{and} \quad
%  \Phi_\bfla(\bfbeta) = \| (D \bfbeta)_+\|_0 .
% \]
% This also belongs to the problem of {\em multiple change-point detection} listed as
% the first applications in \cite{fan2018approximate}. 

%%%%%%%%%%%%%%%%%

%%------------------------------------
\subsection{Research questions and literature review}

The concept of {\em stationary duality} was first proposed in \cite{zhang2025composite}, where the two-block
case (i.e., $g=0$) of \eqref{P} was considered under the 
assumption that $f$ is strongly convex and $\Phi_\bfla = \Phi_{+,\bfla}$ with $\bfla = \lambda \bfone_r$ and $\lambda>0$. This means that 
its conjugate function $f^*$ is continuously differentiable
\cite{rockafellar1970convex}.
We recognize that strong convexity has its specific domain of applications, but it is violated by 
some examples listed above when they are cast as a two-block structural problem.
Relaxation of the strong convexity to convexity immediately raises several issues including
the existence of global solutions for the primal and dual problems. 

%Although this assumption simplifies its analysis, there are extensive applications in forms of three-block problem \eqref{P}. This paper develops its stationary dual problem \eqref{D} based on the stationary duality in \cite{zhang2025composite}. Then several questions naturally arises.
%our contribution: stationary duality of P
%
%
%
%Moreover, several questions were left unexplored.

There is a substantial body of literature on local optimality conditions for CCOP or more general models, typically formulated in terms of suitable stationarity concepts (see, e.g., \cite{li2015global,bot2019proximal,zhang2024zero,bolte2018nonconvex,zhou2021quadratic,cui2023minimization,hananalysis}). However, studies on the existence of global solutions is relatively limited. The first question we address is the {\bf existence of global solutions} for 
both the primal and dual CCOP.
For the special case of \eqref{Energy-minimization} with $\Omega = \Re^n$, and $\Phi_\bfla = \Phi_{0,\bfla}$ (the potential function of ($f11$) in \cite[Table~1.1]{nikolova2005analysis}) and $\bfb=0$, Nikolova characterized its global solution by equivalent convex constrained reformulation under the assumptions that $\Ker A \cap \Ker B = \{0\}$ holds and $B$ is a finite-difference linear operator. 
In another study \cite{nikolova2013description},
Nikolova used the two important tools from optimization \cite{auslender2003asymptotic} to study the existence of global 
solutions of a simpler problem in the setting of $f=0$, $g(\bfu) = \| \bfu-\bfa\|^2/2$, $B = I$, and $\bfb=0$. 
The tools are about a function being {\em Asymptotically Level-set Stable} (ALS) and its asymptotic function
being nonnegative. These tools are also applied to research existence of global solutions of two-block case with $g$ taken as Huber loss \cite{akkaya2020minimizers} and $\ell_1$ loss \cite{akkaya2025minimizers}. Our current work extend them to the three-block case with convex $f$, $g$ (possibly nonsmooth and allowing $\dom f \neq \mbR^n$, and $\dom g \neq \mbR^m$). 
We will also see that both the matrix $B$ and the vector $\bfb$ play important roles in the solution existence.
At a technical level, the results in \cite{akkaya2020minimizers,akkaya2025minimizers,nikolova2013description} only cover the two-block
case with $\Phi_\bfla = \Phi_{0,\bfla}$ and $\bfb =0$. Furthermore, 
they cannot be trivially extended to $\Phi_\bfla = \Phi_{+,\bfla}$ or $\bfb \not=0$. For instance, the problem \eqref{hSVM} is not covered due to $\Phi_\bfla = \Phi_{+,\bfla}$ and $\bfb = -\bfone_r$. Therefore, the existence of its global solution 
is not guaranteed through the existing results. Moreover, we need to stress that all the three works \cite{akkaya2020minimizers,akkaya2025minimizers,nikolova2013description} focus on the existence of global solution to their respective primal problem. The tools developed in these papers can also be applied to the associated stationary dual problems. However, if the relationship between primal and dual problems is further explored, we may obtain more concise sufficient condition for the existence of solutions to the stationary dual problem. For example, in the classic convex optimization theory, the Slater condition of the primal problem is sufficient for the existence of the dual solution (see, e.g. \cite[Proposition 6.4.4]{bertsekas2003convex}), and this condition is often easy to verify. This inspires us to study the existence of global solution to Problem \eqref{D} by exploiting primal-dual relationship.

% Although the developed tools can also be applied to the associated dual problems, the implicit relation between primal and dual problems  

%Moreover, we can handle nonsmooth convex loss functions such as
%$\ell_1$-loss (in comparison to the least-squares).

%With the variable substitution $\bfu := B\bfx - \bfb$, the 
%composite counting is replaced by a simple counting on $\bfu$, but at a price of additional linear constraint: $B \bfx - \bfu = \bfb$. Under this setting, a large number of algorithms for
%nonconvex and nonsmooth optimization can be applied, see e.g.,
%\cite{li2015global, bot2019proximal, zhang2024zero}. Unfortunately, the variable substitution does
%not help in the study of existence of global solutions of the primal problem \eqref{P}.

The second question we need to address is the {\bf solution correspondence} between the primal \eqref{P} and the dual 
\eqref{D}.
Under the two-block case where $f$ is strongly convex, $g=0$, and $\Phi_\bfla = \Phi_{+,\bfla}$, one-to-one correspondence between local solutions of the primal
and the dual problems were established in \cite{zhang2025composite}.
We extend this correspondence to the three-block case with $\Phi_\bfla \in \{ \Phi_{0,\bfla}, \Phi_{+,\bfla} \}$ and only assume the convexity of $f$ and $g$. Furthermore, the correspondence between global solutions of \eqref{P} and \eqref{D} is a challenging question which has not been addressed in \cite{zhang2025composite}.
Positive answer to this question is the foundation for studying
the dual problem, which is a class of cardinality optimization well
surveyed in \cite{tillmann2024cardinality}. This opens
a large venue for research on the composite cardinality optimization.

%%%%%%%
\subsection{Main contributions}

We summarize the main contributions on the three-block model \eqref{P} and its dual problem \eqref{D}.

\begin{itemize}
	\item[(i)] \textbf{On existence of global solutions}. For Problem \eqref{P}, the main technique to ensure the existence of a global solution is the ALS property of the objective function. On one hand, we study the ALS property of the sum of an ALS function and the cardinality term $\Phi_\bfla(B\bfx - \bfb)$, see Lem. \ref{Lemma-ALS}. On the other hand, we show that the convex composition term $g(A\bfx)$ is also ALS, see Lem. \ref{lem-ALS-composition}. In particular, we consider the case $f(\bfx) = \delta(\bfx | C)$, where $C$ is {\em polyhedral convex}, and give a sufficient condition to ensure the existence of a global solution to \eqref{P}, see Thm. \ref{Theorem-Primal}. This case covers the examples introduced in Section \ref{sec-examples} and the models investigated in \cite{nikolova2013description,akkaya2020minimizers,akkaya2025minimizers}. 
	For Problem \eqref{D}, we reformulate it as the minimization of finitely many convex optimization problems, each restricted to a subspace, see Lems. \ref{lem-piecewise-global-phi+} and \ref{lem-piecewise-global-phi0}. With the aid of strong duality theorem and Slater conditions (see, Assumptions \ref{assumption-slater+} and \ref{assumption-slater}) in convex programming, we prove the existence of global solutions to \eqref{D}. 

	\item[(ii)] {\bf On illustration}. In Section~\ref{Subsection-Illustration}, the sufficient conditions for existence of global solutions are applied to the examples described above. 
	This serves three purposes. 
	Firstly, we show these conditions are easily verified for practical applications. Secondly, compared with existing works, our conditions are applicable to more general cases, e.g. nonsmooth $f,g$, 
	allowing $\dom f \neq \mbR^n$,
	$\dom g \neq \mbR^m$, $\Phi_\bfla = \Phi_{+,\bfla}$, and $\bfb \neq 0$.  
	%we also show that we can cover nonsmooth loss functions
	%(e..g, $\ell_1$-loss vs least-squared loss).
	%Section~\ref{Subsection-Illustration} is devoted for the illustration.
	Thirdly, we demonstrate that the proposed sufficient conditions are also necessary for the stationary dual problem \eqref{hSVM}. 
	This suggests that the sufficient conditions are difficult to improve in general. 
	This result provide a theoretical justification for including a ``regularization'' term in the sparse dual SVM proposed by \cite{zhou2021sparse}, where only empirical motivation was given.

	\item[(iii)] {\bf On the one-one correspondence of solutions}. We first characterize the local solutions of \eqref{P} and \eqref{D} by their stationary points. Then through the use of stationary duality, we prove that $\bfx^*$ is local solution of the primal problem if and only if there exists a local solution $(\bfy^*, \bfz^*)$ for the dual problem, see Thm. \ref{thm-correspondence}. The correspondence of global solutions between \eqref{P} and \eqref{D} is a more challenging question. We demonstrate that the correspondence holds for a pair of primal and dual global solutions if the regularization parameter $\bfmu$ of \eqref{D} is appropriately selected, see Remark \ref{remark-correspondence-global}. 
	
	% We also characterize their relationship how
	%they can be derived from each other, see Thm???
	%This result lays the foundation for studying the primal problem
	%via its dual problem.	
\end{itemize}

\subsection{Organization}
In the next section, we collect the tools to be used in our analysis. The standard references are \cite{robinson1976implicit,auslender2003asymptotic,mordukhovich2018variational}. 
We also include the detailed derivation of the dual problem.
Section~\ref{Section-Global-Solution} studies the sufficient
conditions for the existence of global solutions for both
the primal and the dual problems. We also demonstrate how those
conditions are validated through examples. 
In Section~\ref{Section-Solution-Correspondence}, 
we establish one-to-one solution correspondence between the primal and dual problems.
We conclude the paper in Section~\ref{Section-Conclusion}.

%%%%%%%%%%%%%%%%%%%%%%%%%%%%%%%%%%%%%%%%%%%%%%%%%%%%%%
\section{Preliminaries} \label{Section-Preliminaries}

This section aims to describe the notations and definitions frequently used in the paper. 

\subsection{Notation} \label{Subsection-Notation}

We denote $\mbR^n$ as $n$-dimensional Euclidean space endowed with the standard inner product $\langle \cdot, \cdot \rangle$ and the induced norm $\| \cdot \|$ is the $\ell_2$-norm (Euclidean norm). $\Re^n_+$ (resp.~$\Re^n_-$) denotes the nonnegative (resp. nonpositive) orthant. The boldfaced lowercase letter $\bfx \in \mathbb{R}^n$ denotes a column vector of size $n$ and $\bfx^\top$ is its transpose. Let $x_i$ or $\bfx[i]$ denote the $i$-th element of $\bfx$. 
For $\bfx \in \Re^n$, we define its sign vector $\mbox{sgn}(\bfx)$
by $(\mbox{sgn}(\bfx))[i]= 1$ for $\bfx[i] > 0$, $-1$ for
$\bfx[i] < 0$, and $0$ for $\bfx[i] = 0$.
%Every positive semidefinite matrix $H$ induces a seminorm $\| \bfx \|_H:= \sqrt{\langle \bfx, H \bfx \rangle}$, where ``$:=$'' means ``define''. 
The $\epsilon$-neighborhood of $\bfx^* \in \mathbb{R}^n$ is denoted as $\mathcal{N}(\bfx^*, \epsilon) := \{ \bfx \in \mathbb{R}^n \ | \ \| \bfx - \bfx^* \| \le \epsilon \}$. 
Given a set $\Omega$, $\mbox{int}(\Omega)$ and $\ri(\Omega)$ represent the interior and the relative interior of $\Omega$ respectively.
For a given matrix $A$, $\Ker A$ (resp.~$\Rge A$) denotes
its kernel (resp.~range) space. 
For a subset $\U \subset \Re^n$, $\conv(\U)$ is the convex hull of 
$\U$.
%we denote $\dist_H(\bfx,\Omega):= \inf_{\bfy \in \Omega}\| \bfx - \bfy \|_H $ and $\dist(\bfx,\Omega):= \inf_{\bfy \in \Omega}\| \bfx - \bfy \| $. The indicator function $\bfdt_\Omega(\bfx)$ takes 0 when $\bfx \in \Omega$ and $\infty$ otherwise. We use $\ri(\Omega)$ to represent the relative interior of $\Omega$.

Let $[r]$ be the set of indices $\{1, \ldots, r\}$.
For a subset $T \subset [r]$, $|T|$ denotes the number of elements of $T$, and $\oT$ consists of those indices of $[r]$ not in $T$. For vector $\bfz \in \mathbb{R}^r$ (resp. matrix $B \in \mbR^{r \times n}$), $\bfz_{T}$ (resp.~$B_{T:}$) denotes the subvector of $\bfz$ indexed by $T$ (resp. the submatrix of $B$ with rows indexed by $T$). The support set of $\bfz \in \mathbb{R}^r$ is denoted by $\mathcal{I}(\bfz):= \{ i\in [r]: z_i \neq 0  \}$. We denote $I$ as the identity matrix of appropriate dimension. For vectors $\bfu, \bfz \in \mbR^r$, $\bfu \perp \bfz$ means $\langle \bfu, \bfz \rangle = 0$.

%Given a function $\varphi: \mbR^n \to (-\infty, \infty]$, we denote $\dom \varphi:=\{ \bfx ~|~ \varphi(\bfx) < \infty \}$. 
%For a set $\Omega \subset \mbR^m$, ${\rm conv}(\Omega)$ is the convex hull of
%$\Omega$.
%For a symmetric matrix $H \in \mbR^{m \times m}$, $H_{\Gamma}$ is the submatrix with both rows and columns indexed by $\Gamma$.
%Particularly, for a differentiable function $h$, we denote $\nabla_{T} h(\bfz):= (\nabla h(\bfz))_{\Gamma}$.
%The norm $\| \cdot\|$ denotes the $\ell_2$ norm of a vector or matrix. For a positive semidefinite matrix $M$, 

%

\subsection{Subdifferentials and conjugacy}

%The main references for this part are \cite{mordukhovich2018variational,RockWets98} 
Let $\psi: \mbR^n \to (- \infty, \infty]$ be a proper and lower semicontinuous (lsc) function. We denote $\dom \psi:=\{ \bfx ~|~ \psi(\bfx) < \infty \}$. The Fr\'{e}chet subdifferential of $\psi$ at $\obx \in \dom \psi$ is defined by
\begin{align*}
	\whp \psi(\obx):= \left\{ \bfq ~ \left| ~ \liminf_{ \bfx \neq \obx,~ \bfx \to \obx } \frac{\psi(\bfx) - \psi(\obx) - \langle \bfq, \bfx - \obx \rangle }{\| \bfx - \obx \|} \geq 0 \right. \right\}.
\end{align*}
If $\obx \notin \dom \psi$, then $\whp \psi(\obx) = \emptyset$. The limiting (or Modukhovich) subdifferential is defined as follows:
\begin{align*}
	\partial \psi(\obx) := \left\{ \bfq ~\left|~ \exists \bfx^\nu \to \obx, \psi(\bfx^\nu) \to \psi(\obx)~\mbox{and} ~ \psi(\bfx^\nu) \ni \bfq^\nu \to \bfq \right. \right\}.
\end{align*}
When $\psi$ is convex, the definition coincides with the one in convex analysis \cite{robinson1976implicit}.
%\begin{align*}
%	\partial \psi(\obx) = \{ \bfq ~|~ \langle \bfq, \bfx - \obx \rangle \leq \psi(\bfx) - \psi(\obx),~ \forall \bfx \in \mbR^n  \}.
%\end{align*}
By \cite{hu2025iterative,le2013generalized} and the separable structure of $\Phi_{0,\bfla}$, the limiting subdifferential of $\Phi_{0,\bfla}$ can be represented as follows:
\begin{align} \label{subdiff-l0}
	\partial \Phi_{0,\bfla} (\bfu) = \left\{ \bfz ~\left|~ z_i \in \left\{ \begin{aligned}
		&  \mbR, && \mbox{if}~ u_i = 0, \\
		& 0,     && \mbox{otherwise},
	\end{aligned}  ~ i \in [r]  \right.  \right. \right\} 
\end{align}
Noticing that \eqref{subdiff-l0} is independent of parameter $\bfla$, $\partial\Psi_{0,\bfmu}$ has the same representation as \eqref{subdiff-l0}. Through the use of above representation, the following relation holds for any given $\bfla,\bfmu > 0$:
\begin{align} \label{stationary-duality}
	\bfz \in \partial  \Phi_{0,\bfla}(\bfu) ~\iff~ \bfu \in \partial \Psi_{0,\bfmu}(
	\bfz),~\forall~ \bfu,\bfz \in \mbR^r.
\end{align}
According to \cite[Section 2]{zhang2025composite}, limiting subdifferential of $\Phi_{+,\bfla}$ and $\Psi_{+,\bfmu}$ are as follows:
\begin{align} \label{subdiff-phi+}
	\partial \Phi_{+,\bfla} (\bfu) = \left\{ \bfz ~\left|~ z_i \in \left\{ \begin{aligned}
		&  \mbR_+, && \mbox{if}~ u_i = 0, \\
		& 0,     && \mbox{otherwise},
	\end{aligned}  ~ i \in [r]  \right.  \right. \right\} 
\end{align}
\begin{align} \label{subdiff-psi+}
	\partial \Psi_{+,\bfmu} (\bfz) = \left\{ \bfu ~\left|~ u_i \in \left\{ \begin{aligned}
		&  \mbR, && \mbox{if}~ z_i = 0, \\
		& 0,     && \mbox{otherwise},
	\end{aligned}  ~ i \in [r]  \right.  \right. \right\} 
\end{align}
and similar relation holds for them:
\begin{align} \label{stationary-duality+}
	\bfz \in \partial  \Phi_{+,\bfla}(\bfu) ~\iff~ \bfu \in \partial \Psi_{+,\bfmu}(
	\bfz),~\forall~ \bfu,\bfz \in \mbR^r.
\end{align}
Given $\psi:\mbR^n \to (-\infty,\infty]$, the conjugate function \cite{mordukhovich2022convex,rockafellar1970convex} of $\psi$ is defined by
\begin{align*}  %\label{conjugate-fun}
	\psi^*( \bfq ) := \sup_{\bfx \in \mbR^n} \langle \bfq, \bfx \rangle - \psi(\bfx).
\end{align*}
%One can verify the following Fenchel's inequality \cite[Theorem 4.6]{beck2017first} holds for any $\bfx, \bfq \in \mbR^n$:
%\begin{align} \label{fenchel-ineq}
%	\psi(\bfx) + \psi^*(\bfq) \geq \langle \bfx, \bfq \rangle.
%\end{align}
If $\psi$ is proper, lsc, and convex, the following claims \cite[Thm.~23.5]{rockafellar1970convex} are equivalent for $\bfx,\bfq \in \mbR^n$:
\begin{align} \label{conjugate-subgradient}
	\psi(\bfx) + \psi^*(\bfq) = \langle \bfx, \bfq \rangle ~\iff~ \bfq \in \partial \psi(\bfx) ~\iff~ \bfx \in \partial \psi^*(\bfq).
\end{align}
The first equation above is often called Fenchel-Young equality.

%Furthermore, given $\sigma > 0$, the conjugacy for $\sigma$-strong convexity and $(1/\sigma)$-smoothness (\cite[Definitions 5.1 and 5.16]{beck2013sparsity}) has the following correspondence \cite[Theorem 5.26]{beck2013sparsity}:
%\begin{align*}
%	\sigma\mbox{-strong convexity of}~\psi \iff (1/\sigma)\mbox{-smoothness of}~\psi^*.
%\end{align*}

%Since the cardinality function is nonconvex and bounded, its conjugate function may not be helpful in theory or computation. Actually, we can verify that $(\mu \card{\cdot})^* = \bfdt_{\{0\}}(\cdot)$, which only be finite at origin. However, one can observe from \eqref{stationary-duality} that the cardinality function and itself satisfies the second relationship in \eqref{conjugate-subgradient}. This will be useful to establish the dual formulation of \eqref{P} and design efficient algorithms. 

%%%%-------------------------------------------
\subsection{Asymptotic functions}

For  a given proper function $\psi: \Re^n \to (-\infty, \infty]$, which is also lsc, its asymptotic function and the associated kernel are given by
\[
\psi_\infty (\bfd) := \liminf_{\stackrel{\widetilde{\bfd} \rightarrow \bfd }{t \rightarrow + \infty} } \frac{\psi ( t \widetilde{\bfd})}{t} \quad \mbox{and} \quad
\Ker \psi_\infty := \left\{ \bfd \ | \ \psi_\infty (\bfd) = 0\right\}.
\]
For a given set $C \subset \Re^n$, its asymptotic cone is defined as
\[
C_\infty := \left\{
\bfd \in \Re^n  \left| 
\mbox{there are sequences} \ t_k \rightarrow + \infty, \
\bfx^k \in C \ \mbox{such that} \
\lim_{k \rightarrow \infty} \frac{\bfx^k}{t_k} = \bfd 
\right.  \right\}.
\]
The definitions were well studied in \cite[Chapter~2]{auslender2003asymptotic}.
When $\psi$ and $C$ are convex, the concepts become the more familiar recession function and recession cone \cite{rockafellar1970convex}.

The following results are collected from \cite{auslender2003asymptotic}.

\begin{lemma} \label{Lemma-Recession-Function}
	The following results hold.
	\begin{itemize}
		\item[(i)] \cite[Cor.~2.5.1]{auslender2003asymptotic}
		For a nonempty set $C \subseteq \Re^n$ one has
		$(\delta_C)_\infty   = \delta_{C_\infty}$.
		
		%	 \cite[Thm.~2.5.4]{auslender2003asymptotic}
		%	Let $f: \Re^n \to (-\infty, \infty])$ be proper and convex. Then $(f^*)_\infty = \sigma_{\dom f}$.}
	
	\item[(ii)] \cite[Prop.~2.6.1]{auslender2003asymptotic}
	Let $f_1, f_2: \Re^n \to (-\infty, \infty]$ be lsc, proper and convex with $\dom f_1 \cap \dom f_2 \not=\emptyset$. 
	Then $(f_1+f_2)_\infty = (f_1)_\infty + (f_2)_\infty$.
	
	\item[(iii)] \cite[Prop.~2.6.3]{auslender2003asymptotic}
	Let $g: \Re^n \to (-\infty, \infty]$ be proper, lsc, and convex. Denoting $\overline{g} (\bfx) := g(A\bfx)$, if $\Rge A \cap \dom g \not= \emptyset$.
	Then  $\overline{g}_\infty(\bfd) = g_\infty (A\bfd)$ for $\bfd \in \Re^n$.
	
\end{itemize}
\end{lemma}

\subsection{Deriving the dual problem} \label{Subsection-Deriving}
To derive the dual problem for \eqref{P} with $\Phi_\bfla \in \{ \Phi_{0,\bfla}, \Phi_{+,\bfla}\}$,
we follow the standard way for a general optimization problem through its parameterization \cite[Chp.~11]{RockWets98}. Let
\[
\V(\bfx, \bfu, \bfv) := f(\bfx) + g(A\bfx + \bfu) + \Phi_\bfla(B \bfx - \bfb + \bfv) .
\] 
The Fenchel conjugate of $\V(\bfx, \bfu, \bfv)$ is
\begin{align*}
\V^*(\bfq, \bfy, \bfz) &= \sup_{\bfx, \bfu, \bfv} 
\left\{
\langle \bfx, \bfq \rangle + \langle \bfu, \bfy \rangle
+ \langle \bfz, \bfv \rangle
- f(\bfx) - g(\underbrace{A\bfx + \bfu}_{=: \bfs}) - \Phi_\bfla(\underbrace{B \bfx - \bfb + \bfv}_{:= \bft})
\right\}	\\
&= \sup_{\bfx, \bfs, \bft} 
\left\{
\langle \bfx, \bfq \rangle + \langle \bfs-A\bfx, \bfy \rangle
+ \langle \bft-B\bfx+\bfb, \bfz \rangle
- f(\bfx) - g(\bfs) -  \Phi_\bfla(\bft)
\right\} \\
&=  \langle \bfb, \bfz \rangle +
\underbrace{\sup_{\bfx} \left\{
	\langle \bfq - A^\top \bfy - B^\top \bfz, \bfx \rangle - f(\bfx)
	\right\}}_{= f^*( \bfq - A^\top \bfy - B^\top \bfz)}
+ \underbrace{\sup_{\bfs} \left\{
	\langle \bfy, \bfs \rangle  - g(\bfs)
	\right\}}_{= g^*(\bfy) } \\
& \quad + \underbrace{\sup_{\bft} \left\{
	\langle \bfz, \bft \rangle  - \Phi_\bfla(\bft)
	\right\}}_{\Phi_\bfla^*(\bft)}\\
&= \langle \bfb, \bfz \rangle 
+ f^*( \bfq - A^\top \bfy - B^\top \bfz)
+ g^*(\bfy) 
+ \Phi_\bfla^*(\bfz).
\end{align*}
The classical dual problem is then given by
\begin{equation} \label{fenchel-duality}
\sup_{\bfy, \bfz} - \V^*(0, \bfy, \bfz)
=  \sup_{\bfy, \bfz} \left\{
- \langle \bfb, \bfz \rangle 
- f^*(- A^\top \bfy - B^\top \bfz)
- g^*(\bfy) 
- \Phi_\bfla^*(\bfz)
\right\} .
\end{equation}
However, due to the boundedness of $\Phi_\bfla(\cdot)$, its conjugate $\Phi_\bfla^*(\bfz)$ is
$0$ when $\bfz=0$ and $+\infty$ otherwise. 
This renders the classical dual problem completely invalid
without inheriting any useful information from $\Phi_\bfla(\cdot)$. Therefore, the direct extension of Fenchel conjugacy for $\Phi_\bfla$ is not suitable.
An alternative dual formulation was recently proposed in \cite{zhang2025composite}.
The main idea is motivated by the conjugate subdifferential correspondence in \eqref{conjugate-subgradient}. Specifically, we replace $\Phi_\bfla^*$ with a function $\Psi_\bfmu(\cdot)$ satisfying
the stationary duality relationship:
\begin{equation}
\label{stationary-duality-phi-psi}
\bfz \in \partial \Phi_\bfla(\bfu) \quad 
\mbox{if and only if} \quad
\bfu \in \partial \Psi_\bfmu(\bfz) \qquad \forall\ \bfu, \bfz \in \Re^r.
\end{equation}
%Therefore, the dual problem of \eqref{P} is then given by (in the form of minimization):
%\begin{align}  \tag{D} \label{D}
%	\min_{(\bfy, \bfz) \in \mbR^m \times \mbR^r} \
%	G(\bfy, \bfz) := \langle \bfb, \bfz \rangle 
%	+ f^*(-A^\top \bfy -B^\top \bfz) + g^*(\bfy) + \mu \Psi_\bfmu(\bfz),
%\end{align}
%where $\mu>0$. It is easy to see that the corresponding stationary dual functions for $\Phi_{0,\bfla}$ and $\Phi_{+,\bfla}$ are respectively given by
%\[
%  \Psi_{0,\bfmu}(\bfz) = \| \bfz\|_0 \qquad \mbox{and} \qquad 
%  \Psi_{+,\bfmu}(\bfz) = \delta_{\Re^r_+} (\bfz) \quad \mbox{(indicator function)}.
%\]
%In terms of \eqref{conjugate-subgradient}, we dropped the requirement of the Fenchel-Young equality and we only kept the
%conjugate subdifferential correspondence.
We call such $\Psi_\bfmu$ a {\em stationary dual} function of $\Phi_\bfla$. As we demonstrate in \eqref{stationary-duality} and \eqref{stationary-duality+}, $\Psi_\bfmu\in \{ \Psi_{0,\bfmu}, \Psi_{+,\bfmu}\}$ and $\Phi_\bfla \in \{ \Phi_{0,\bfla}, \Phi_{+,\bfla}\}$ satisfy the correspondence \eqref{stationary-duality-phi-psi}.

Heuristically, we replace $\Phi_\bfla^*$ in \eqref{fenchel-duality} by the associated stationary dual function $\Psi_\bfmu \in \{ \Psi_{0,\bfmu}, \Psi_{+,\bfmu} \}$. In this way, we obtain our dual problem \eqref{D}, which is presented in the form of minimization. The theoretical justifications for this dual formulation are provided in Section \ref{Section-Solution-Correspondence}. Specifically, we show that $\bfx^*$ is a local solution of \eqref{P} if and only if there exists a local solution $(\bfy^*, \bfz^*)$ for \eqref{D}, which means the correspondence of primal and dual local solutions. If the weighted parameter $\bfmu$ in $\Psi_\bfmu$ is appropriately selected, we further demonstrate that the correspondence holds for a pair of global solutions to the primal and dual problems.

%  and we used the relationship 
%$\partial (\mu \Psi_\bfmu) = \partial \Psi_\bfmu$ for any $\mu>0$. 
%This allows us to introduce $\mu$ in the dual problem.
%The rest of the paper is to justify why the dual problem provides a viable approach to
%the primal problem through addressing the two key research
%questions discussed in Introduction.

%%%%%%%%%%%%%%%%%%%%%%%%%%%%%%%%%%%%%%%%%
\section{Existence of Global Solutions} \label{Section-Global-Solution}

\subsection{A general existence result}

We first state a general existence result of global solutions for the optimization 
problem:
\be \label{General-Problem}
\inf_{\bfx \in \Re^n} \ \cL(\bfx),
\ee 
where $\cL: \Re^n \to (-\infty, \infty]$ is an lsc, proper function.
We define $S_\cL$ to be the set of sequences of the following property:
\[
S_\cL := \left\{
\{ \bfx^k\} \left| \| \bfx^k\| \rightarrow \infty, \
\{ \cL(\bfx^k)\} \ \mbox{bounded above and} \
\bfx^k /\| \bfx^k \| \rightarrow \overline{\bfx} \in \Ker \cL_\infty
\right. \right\}.
\]
Certainly, when $\cL(\cdot)$ is level-bounded $(\cL(\bfx) \rightarrow +\infty$ as
$\| \bfx \| \rightarrow \infty$), $S_\cL = \emptyset$.
When $S_\cL$ is not empty, we need $\cL(\cdot)$ to be {\em asymptotically level stable} (ALS). This definition involves the level set of $\cL$ at a level $\tau \in \Re$:
\[
\mbox{lev}(\cL, \tau) := \left\{
\bfx \in \Re^n \ | \ \cL(\bfx) \le \tau
\right\} .
\]

\begin{definition} \label{Def-ALS}
\cite[Def.~3.3.2]{auslender2003asymptotic}
Let $\cL: \Re^n \to (-\infty, \infty]$ be an lsc, proper function.
Then $\cL$ is said ALS if for each $\rho>0$, each bounded sequence of reals ${\tau_k}$ and each sequence $\{ \bfx^k\}$ satisfying
\[
\bfx^k \in \mbox{\rm lev}(\cL, \tau_k), \quad
\| \bfx^k \| \rightarrow \infty, \quad
\bfx^k / \| \bfx^k\| \rightarrow \overline{\bfx} \in \Ker \cL_\infty,
\]
there exists $k_0$ such that $\bfx^k - \rho \overline{\bfx} \in \lev (\cL, \tau_k)$ for all $k \ge k_0$.
\end{definition}
We now give two conditions ensuring existence of a global solution to Problem
\eqref{General-Problem}:
\begin{flalign*}
&\text{(C0)} \quad \cL_\infty (\bfd) \ge 0, \ \forall \ \bfd \in \Re^n \backslash \{0\}. &\\
&\text{(C1)} \quad  
\mbox{Either} \ S_\cL = \emptyset \ \mbox{or} \
\cL \ \mbox{is ALS}. &
\end{flalign*}
%\[ %\lefttag{$\H_0$}
% 
%\]
%
%\[
%
%\]
%We now apply the result \cite[Cor.~3.4.2]{auslender2003asymptotic} to get the following
%result. 
We need the following result, which is key to our derivation on existence of global solutions.

\begin{theorem} \label{Thm-H01}
\cite[Cor.~3.4.2]{auslender2003asymptotic}
Suppose the two conditions (C0) and (C1) hold. Then
\eqref{General-Problem} has a global solution.	
\end{theorem}

%\begin{remark}\label{Remark-Assumption-H1}
In fact, there is Condition $\H_1$ proposed in \cite{auslender2003asymptotic}, which together with (C0), forms a necessary and sufficient condition for \eqref{General-Problem} to have a global solution. (C1) being part of the sufficient condition 
is easier to verify and it suits our purpose.
%For example, $\cL(\cdot)$ being weakly coercive is sufficient for
%the existence of a global solution.	
%\end{remark}
%%%%%%%%%%%%%%%%%%%%%%%%%%%%%%%%%%%%%%%%%%%%%%%%%%%%%%%%%
\subsection{Global solution of the primal problem}

%For a special two-block case
%\[
%   \min_\bfx \frac 12 \| A\bfx - \bfa \|^2 + \lambda \| \bfx\|_0 ,
%\]
%where $\lambda > 0$, a global solution exists \cite[Sect.~4]{nikolova2011energy}.
%We will see that this result can not be easily extended to three-block case.

The following technical results may be independently interesting to the study of ALS functions.

\begin{lemma} \label{Lemma-ALS}
Let $\cL: \Re^n \to (-\infty, \infty]$ be an lsc and proper function.
The following assertions hold.

\begin{itemize}
	\item[(i)] If $\cL(\cdot)$ is ALS, then so is $H(\bfx) := \cL(\bfx) + \Phi_\bfla( B\bfx - \bfb)$, where $\Phi_\bfla \in \{ \Phi_{0,\bfla}, \Phi_{+,\bfla}\}$. 
	
	\item[(ii)] If $\cL(\cdot)$ is ALS and convex, then so is $U(\bfx) := \cL(\bfx) + \delta(\bfx | C)$, provided that
	$C$ is a polyhedral convex set and
	$\dom \cL \cap C \not= \emptyset$.
	%  and   $\cL_\infty(\bfd) > - \infty$ for any $\bfd \not\in C_\infty$.
\end{itemize}
\end{lemma} 

\bp
(i) 
We first note that by definition
\[ \label{asym-H}
H_\infty (\bfd) = \liminf_{\stackrel{\widetilde{\bfd} \rightarrow \bfd}{ t\rightarrow + \infty}} \frac{ H(t \widetilde{\bfd} )}{t} = \cL_\infty (\bfd)
\]
because $\Phi_\bfla$ is bounded (i.e., $0\le \Phi_\bfla\le \sum_{i=1}^r \lambda_i$). Therefore, 
$
\Ker H_\infty = \Ker \cL_\infty .
$
For the proof, we suppose $S_H$ is not empty and $\rho>0$ is given. 
Let $\{(\bfx^k, \tau_k)\}$ be a 
sequence such that $\{ \tau_k\}$ is bounded above and
\be \label{S-Sequence}
\bfx^k \in \lev(H, \tau_k), \quad
\| \bfx^k\| \rightarrow \infty, \quad
\bfx^k/ \| \bfx^k\| \rightarrow \overline{\bfx} \in \Ker H_\infty.   %= \Ker \cL_\infty.
\ee
Define $\beta_k := \tau_k -  \Phi_\bfla( B\bfx^k - \bfb)$.
Then we have
\[
\cL(\bfx^k) = H(\bfx^k) - \Phi_\bfla( B\bfx^k - \bfb)
\le \tau_k - \Phi_\bfla( B\bfx^k - \bfb) = \beta_k.
\]
That is,  $\bfx^k \in \lev( \cL, \beta_k)$ and $\{ \beta_k\}$ is also bounded above due to the uniform boundedness of $\Phi_\bfla(\cdot)$.
Since $\cL$ is ALS, we must have $\bfx^k - \rho \overline{\bfx} \in \lev(\cL, \beta_k)$ for all $k \ge k_0$ for some $k_0$. here we used the fact $\overline{\bfx} \in \Ker H_\infty= \Ker \cL_\infty$.
Now we consider the sequence
$\{ \bfx^k - \rho \overline{\bfx} \}$ and consider two cases depending on the value of 
$\overline{\bfu} := B \overline{\bfx}$. 

{\bf Case 1.} $\overline{\bfu}[i] \not=0$. Using the limit $\lim_{k \rightarrow \infty} B\bfx^k/\| \bfx^k\| = B \overline{\bfx} = \overline{\bfu}$, we have $| (B \bfx^k)[i] | \rightarrow \infty$.
Consequently, we must have the following results for $k$ sufficiently large
\[
(  B( \bfx^k - \rho \overline{\bfx} ) - \bfb ) [i] \not= 0 \quad 
\mbox{and} \quad
\mbox{sgn} (  B( \bfx^k - \rho \overline{\bfx} ) - \bfb ) [i]
= \mbox{sgn} ( B \bfx^k - \bfb)[i] .
\]

{\bf Case 2.}
$\overline{\bfu}[i] =0$. We must have $(  B( \bfx^k - \rho \overline{\bfx} ) - \bfb ) [i] = ( B \bfx^k - \bfb)[i]$.

For both cases, we have for $k$ sufficiently large 
$
\Phi_\bfla (B( \bfx^k - \rho \overline{\bfx} ) - \bfb )  
= \Phi_\bfla ( B \bfx^k - \bfb) = \tau_k - \beta_k .
$
Therefore,
\[
H(\bfx^k - \rho \overline{\bfx} ) 
= \cL( \bfx^k - \rho \overline{\bfx}) + \Phi_\bfla (B( \bfx^k - \rho \overline{\bfx} ) - \bfb )
\le \beta_k + (\tau_k -  \beta_k) \le \tau_k,
\]
which implies $\bfx^k - \rho \overline{\bfx} \in \lev(H, \tau_k)$ for $k$ sufficiently large.
This proves that $H$ is ALS. 

(ii) It follows from Lemma~\ref{Lemma-Recession-Function}(i) and (ii) that
\[ \label{af-sum-rule}
U_\infty(\bfd) = \delta_\infty(\bfd | C)  + \cL_\infty (\bfd)
= \delta(\bfd | C_\infty) +  \cL_\infty (\bfd),
\]
The convexity of $\cL$ implies $\cL_\infty (\bfd) > - \infty$ for any $\bfd \in \Re^n$. Therefore, the above equality indicates
\[
\Ker U_\infty = C_\infty \cap \Ker\cL_\infty .
\]
Suppose we have a sequence $\{ (\bfx^k, \tau_k) \}$ satisfying 
\be \notag 
\bfx^k \in \lev(U, \tau_k), \quad
\| \bfx^k\| \rightarrow \infty, \quad
\bfx^k/ \| \bfx^k\| \rightarrow \overline{\bfx} \in \Ker U_\infty.
\ee
Then we have 
\[
U(\bfx^k) = \delta(\bfx^k|C) + \cL(\bfx^k) = \cL(\bfx^k) \le \tau_k,
\]
which means that $\bfx^k \in \lev(\cL, \tau_k)$ and $\bfx^k \in C$. Since $\cL(\cdot)$ is ALS, given $\rho > 0$, we have $\bfx^k - \rho \overline{\bfx} \in \lev(\cL, \tau_k)$ for sufficiently large $k$. Moreover, since $C$ is a polyhedral convex set, it is also asymptotically linear by \cite[Proposition 2.3.1]{auslender2003asymptotic}. Then according to \cite[Definition 2.3.1]{auslender2003asymptotic}, $\bfx^k - \rho \overline{\bfx} \in C$ when $k$ is large enough. Overall, we can obtain $\bfx^k - \rho \overline{\bfx} \in \lev(U, \tau_k)$, which means $U$ is ALS. 
\hfill $\Box$

Next, we point out a subtle difference from a known result.

\begin{remark} \label{Remark-ALS}
Regarding (ii), there is a similar result in \cite[Proposition~3.3.3(c)]{auslender2003asymptotic},
where the same claim holds provided $\cL_\infty (\bfd) \ge 0$ for $\bfd \in C_\infty \backslash \{0\}$. On comparison, we assume that $\cL$ is convex and $C$ is polyhedral convex in this case.
%This latter assumption was implicitly assumed therein. 
%In this case, we use convexity of .
%though it is necessary for $F(\cdot)$ to have a global minimum (i.e., attained at some point). 
Let us use an example to show the improvement.
Let
\[
\cL(\bfx) = \bfx^\top Q \bfx + \langle \bfa, \bfx \rangle
\] 
for some positive semidefinite matrix $Q$ with $\Ker Q \neq 0$ and $\bfa \in \Re^n$. It is ALS according to \cite[Proposition~3.3.3(d)]{auslender2003asymptotic}. 
We also know from \cite[Page~68]{rockafellar1970convex} that
\[
\cL_\infty(\bfd) = \left\{
\begin{array}{ll}
	\langle \bfa, \bfd \rangle & \mbox{if} \ Q\bfd = 0 \\
	+\infty & \mbox{otherwise}.
\end{array} 
\right .
\]
Let $C =\left\{ \bfx \; | \; Q\bfx =0\right\}$. Then, we see that $C_\infty = C$ and
$\cL_\infty (\bfd) = +\infty$ for $\bfd \not\in C_\infty$. 
From Lemma~\ref{Lemma-ALS}(ii), we can claim
that $F(\bfx) = \cL(\bfx) + \delta(\bfx | C)$ is ALS. However, this result cannot be claimed
from  \cite[Proposition 3.3.3(c)]{auslender2003asymptotic} because $\cL_\infty (\bfd) = \langle \bfa, \bfd \rangle$ may be negative when $\bfd \in C_\infty$. Therefore, the claim in (ii)
is an extension of the mentioned result.
\end{remark}

\begin{remark} \label{Remark-Polyhedracity} 
The polyhedracity of $C$ is important. Here is a counterexample in two dimensions. Let $\bfx :=(x_1, x_2)^\top\in \Re^2$ and we denote 
\[
C := \left\{ (x_1, x_2) \ | \ x_2 \ge \exp(x_1)\right\}, \quad 
\cL(\bfx) := x_2^2, \quad
U(\bfx) := \delta(\bfx | C) + \cL(\bfx).
\] 
$C$ is a curved convex region and it is not polyhedral.
Using \eqref{af-sum-rule}, we can compute
\[
\Ker U_\infty = C_\infty \cap \Ker \cL_\infty = \left\{ (t , 0) \ | \ t \leq 0  \right\}.
\]
Let us consider the sequence $\{ \bfx^k\} := \{ (t_k, \exp(t_k)) \}$ with $t_k \rightarrow - \infty$. Then $\bfx^k /\| \bfx^k\| \rightarrow \overline{\bfx} := (-1; 0) \in \Ker U_\infty$ and
$\bfx^k \in \lev(U, \tau_k)$ with $\tau_k = \exp(t_k) < 1$ (due to $t_k <0$).
For any $\rho >0$, consider the points
$\bfx^k - \rho (-1, 0) = (t_k + \rho, \exp(t_k))$. It is easy to see
\[
\exp(t_k + \rho)  = \exp(\rho) \exp(t_k) >  \exp(t_k).
\]
This means that $(t_k + \rho, \exp(t_k)) \not\in C$. 
Hence, $U (\bfx^k - \rho \overline{\bfx}) = +\infty$. That is, $(\bfx^k - \rho \overline{\bfx}) \not \in \lev(U, \tau_k)$. Consequently, $U$ is not ALS even
though $\cL$ is ALS.
\end{remark}

To proceed, we give the following lemma for ALS property of a composite function.
\begin{lemma} \label{lem-ALS-composition}
Let $g: \Re^n \to (-\infty, \infty]$ be proper, lsc, and convex. Denoting $\overline{g} (\bfx) := g(A\bfx)$. If $g$ is ALS, then so is $\overline{g}$.	
\end{lemma}
\bp
Let us consider bounded reals $\tau_k$ and sequence $\{ \bfx^k \}$ satisfying 
\begin{align*}
\bfx^k \in \lev( \overline{g}, \tau_k ), \quad \| \bfx^k \| \to \infty, \quad \bfx^k/\| \bfx^k \| \to \obx \in \Ker \overline{g}_\infty
\end{align*}
On one hand, it holds that $g(A\bfx^k) \leq \tau_k$. On the other hand, we have $A\obx \in \Ker g_\infty$ due to $\overline{g}_\infty(\bfd) = g_\infty(A\bfd) $ from Lem.~\ref{Lemma-Recession-Function}. Given $\rho > 0$, let us consider the following two cases.

{\bf Case 1.} If $A\obx = 0$, then $\overline{g}(\bfx^k - \rho \obx) = g(A\bfx^k) \leq \tau_k$.

{\bf Case 2.} If $A\obx \neq 0$, then we can derive 
\begin{align*}
A\bfx^k \in \lev(g, \tau_k), \quad \| A\bfx^k \| \to \infty, \quad \frac{A\bfx^k}{\|A\bfx^k \|} = \frac{A(\bfx^k/ \| \bfx^k \|)}{\|A(\bfx^k/ \| \bfx^k \|) \|} \to \frac{A\obx}{\| A\obx \|} \in \Ker g_\infty,
\end{align*}
where the last relation holds because $g_\infty$ is homogeneous (see \cite[Proposition 2.5.1]{auslender2003asymptotic}). Denoting $\overline{\rho} := \rho \| A \obx \|$, we can use the ALS property of $g$ to derive $\overline{g} (\bfx^k - \rho \obx) = g(A\bfx^k - \overline{\rho} ( A\obx )/\| A\obx \|) \leq \tau_k$.

The above two cases indicate that $\overline{g}$ is ALS.  
\ep

The main result in this section is an easy consequence.

\begin{theorem} \label{Theorem-Primal}
Let $f (\bfx) = \delta(\bfx|C)$ for some convex set $C \subseteq \Re^n$ and $g: \Re^m \to \Re \cup \{+\infty\}$ be proper, lsc, and convex.
Suppose the following four conditions hold: 
(i) $A C \cap \dom g \not= \emptyset$,
(ii) $g$ is ALS, 
(iii) $C$ is polyhedral, and
(iv) $g_\infty (A\bfd) \ge 0$ for all $\bfd \in C_\infty\backslash \{ 0 \}$.
Then the primal problem \eqref{P} has a global solution.
\end{theorem}

\bp
By the conditions (i)-(iii), and Lems. \ref{Lemma-ALS}(ii) and \ref{lem-ALS-composition}, we get that
$\widetilde{F}(\bfx) := \delta(\bfx| C) + g(A\bfx)$ is ALS. 
By Lem.~\ref{Lemma-ALS}(i), $F(\bfx)=\widetilde{F}(\bfx) + \Phi_\bfla(\bfx)$ is also ALS.
Due to the convexity of both $C$ and $g(\cdot)$, we have
from Lem.~\ref{Lemma-Recession-Function}
\[
F_\infty(\bfd) = \delta(\bfd | C_\infty) + g_\infty(A\bfd) \ge 0, \quad \forall \ \bfd \neq 0 .
\]
%due to the assumption (iii). 
It follows from the general result Thm.~\ref{Thm-H01} that the primal problem has a global solution. 
\hfill $\Box$

Regarding the condition (ii) in the above theorem, it is shown in \cite{auslender2003asymptotic} that piecewise linear quadratic and convex functions, level-bounded functions, and asymptotically linear functions are ALS. The condition (iv) automatically holds for any bounded below functions. Therefore, Thm. \ref{Theorem-Primal} is general enough to cover many practical problems. For example,
$C$ can be a box constraint such as $[L, U]$, $\Re^n_+$, or any polyhedral convex set. 
The function $g(\bfv)$ can take the $\ell_p$ norm $\| \bfv - \bfa\|_p$ ($p\ge 1$), 
$\| \bfv - \bfa\|^2$, Huber loss, or any strongly convex function. This covers the models studied in \cite{akkaya2020minimizers,akkaya2025minimizers,nikolova2013description,nikolova2005analysis} and its extension constrained by a polyhedral
convex set.

%%%%%%%%%%%%%%%%%%%%%%%%%%%%%%%%%%%%%%%%%%%%%
\subsection{Global solution for the dual problem}

In the last section, we mainly use the nonnegativity of the asymptotic function and ALS property to ensure the existence of global solutions of Problem \eqref{P}. By treating \eqref{D} as an independent problem, we can directly apply these tools to verify the existence of dual optimal solution based on the conjugate functions $f^*$ and $g^*$, and other data. 
We choose not to do so because of the following reasons. 
Firstly, sufficient conditions enforced on 
the conjugate functions
$f^*$ and/or $g^*$ are more involved and are often hard to verify.  
Secondly, it is often preferable to use the original data (i.e. functions $f$ and $g$)
to judge whether the dual problem is solvable. 
%However, the primal and dual problems are often closely related. Moreover, we may prefer directly using the original data (i.e. functions $f$ and $g$) in Problem \eqref{P} to judge whether the dual problem is solvable. 
For example, in the classic convex optimization theory, the Slater condition of the primal problem implies the existence of a dual optimal solution (see, e.g. \cite[Proposition 6.4.4]{bertsekas2003convex}). 
Based on the above considerations, we aim to establish the existence of solutions to the dual problem with the aid of primal-dual relation.
%ALS, asymptotic function 
%
%apply to f*, g*
%
%direct use of function data
%
%

Since the cardinality function is constant on finitely many convex regions, we reformulate Problem \eqref{D} as the minimization over a finite collection of convex programs. For each convex program, we can apply the classic strong duality theorem to ensure the existence of a global solution. These form the main procedure of the proof in this section. For convenience, we denote
\begin{align*}
\bfw := [\bfy;\bfz] \quad \mbox{and} \quad \Xi(\bfw) := f^*(-A^\top \bfy - B^\top \bfz) + g^*(\bfy) + \langle \bfb, \bfy \rangle
\end{align*}
%%%%%
\subsubsection{\texorpdfstring{Sufficient conditions for $\Psi_\bfmu = \Psi_{+,\bfmu}$}{Sufficient conditions for Psi_lambda = Psi_+,lambda}}
The subsequent lemma presents an equivalent formulation of Problem \eqref{D}.
\begin{lemma} \label{lem-piecewise-global-phi+}
For Problem \eqref{D} with $\Psi_\bfmu = \Psi_{+,\bfmu}$, it holds that 
\begin{equation} \label{piecewise-global-phi+}
	\inf_{\bfw} \Xi(\bfw) + \Psi_{+,\bfmu}(\bfz) = \inf_{\S \subseteq [r]} \Big\{  \inf_\bfw \Big\{\Xi(\bfw) + \sum_{i \in \S} \mu_i ~\Big|~ \bfz_{\S} \geq 0,~ \bfz_{\omS} = 0 \Big. \Big\} \Big\} 
\end{equation}
\end{lemma}
\bp
\uline{We first prove that ``$\geq$" holds in \eqref{lem-piecewise-global-phi+}.} Given any $\whw = [\why; \whz]$, we have 
\begin{align*}
\Xi(\whw) + \Psi_{+,\bfmu}(\whz) = & \Xi(\whw) + \sum_{i \in \I(\whz)} \mu_i + \delta(\whz | \mbR^m_+) \\
\geq & \inf_\bfw \Big\{\Xi(\bfw) + \sum_{i \in \I(\whz)} \mu_i ~\left|~ \bfz_{\I(\whz)} \geq 0,~ \bfz_{\OI(\whz)} = 0 \right. \Big\} \\
\geq & \inf_{\S \subseteq [r]} \Big\{  \inf_\bfw \Big\{\Xi(\bfw) + \sum_{i \in \S} \mu_i ~\Big|~ \bfz_\S \geq 0,~ \bfz_{\omS} = 0 \Big. \Big\} \Big\} 
\end{align*}
By taking minimization on both side of this inequality with respective to $\whw$, we can derive ``$\geq$" in \eqref{piecewise-global-phi+}.

\uline{Next we prove ``$\leq$" holds in \eqref{piecewise-global-phi+}.} Given $\S \subseteq [r]$ and $\bfw = [\bfy; \bfz]$ with $\bfz_\S \geq 0$ and $\bfz_{\omS} = 0$, the following relations hold:
\begin{align*}
\Xi(\bfw) + \sum_{i \in \S} \mu_i \geq \Xi(\bfw) + \Psi_{+,\bfmu} ( \bfz ) \geq \inf_{\bfw} \Xi(\bfw) + \Psi_{+,\bfmu}(\bfz).
\end{align*}
Minimization on both side of the inequality over $\bfw = [\bfy; \bfz]$ subject to $\bfz_\S \geq 0$ and $\bfz_{\omS} = 0$, and then over $\S \subseteq [r]$, we can derive the desired conclusion.
\ep

To ensure the existence of a solution to Problem \eqref{D}, it suffices to show that for each $\S \subseteq [r]$, the convex program on the right-hand side of \eqref{lem-piecewise-global-phi0} admits a global solution. For this, we need the following assumption.
\begin{assumption} \label{assumption-slater+}
The objective function $F$ of Problem \eqref{P} is bounded below. Moreover, it holds that $\{ \bfx ~|~ \bfx \in \ri(\dom f),~ A\bfx \in \ri(\dom g),~ B\bfx \leq \bfb \} \neq \emptyset$. Particularly, ``\ri" in the set can be omitted if the corresponding function $f$ or $g$ is polyhedral.
\end{assumption}

\begin{theorem} \label{thm-global-sol-D+}
Given $\Psi_\bfmu = \Psi_{+,\bfmu}$, if Assumption \ref{assumption-slater+} holds, then Problem \eqref{D} has a global solution.
\end{theorem}
\bp
The Slater condition of the following convex program holds due to Assumption \ref{assumption-slater+}:
\begin{align} \label{P-S}
\inf_{\bfx} f(\bfx) + g(A\bfx) ~s.t.~ (B\bfx - \bfb)_\S \leq 0.
\end{align}
The optimal value of the above problem must be finite because of the lower boundedness of $F$. The dual problem of \eqref{P-S} can be represented as
\begin{align*}
\inf_\bfw \Xi(\bfw) ~s.t.~ \bfz_\S \geq 0,~\bfz_{\omS} = 0. 
\end{align*}
It follows from \cite[Proposition 6.4.4]{bertsekas2003convex} that the dual problem must have global solutions. Finally, by using Lem. \ref{lem-piecewise-global-phi+}, we can conclude that \eqref{D} must have a global minimizer.
\ep

\subsubsection{\texorpdfstring{Sufficient conditions for $\Psi_\bfmu = \Psi_{0,\bfmu}$}{Sufficient conditions for Psi_bfla = Psi0,bfla}} \label{sec-sufficient-condition0}
The procedure is analogous to that in Section \ref{sec-sufficient-condition0}.
\begin{lemma} \label{lem-piecewise-global-phi0}
For Problem \eqref{D} with $\Psi_\bfmu = \Psi_{0,\bfmu}$, it holds that 
\begin{equation} \label{piecewise-global}
	\inf_{\bfw} \Xi(\bfw) + \Psi_{0,\bfmu}(\bfz) = \inf_{\S \subseteq [r]} \Big\{  \inf \Big\{\Xi(\bfw) + \sum_{i \in \S} \mu_i ~\Big|~ \bfz_{\omS} = 0 \Big. \Big\} \Big\} 
\end{equation}
\end{lemma}
The proof of Lem. \ref{lem-piecewise-global-phi0} is similar to that of Lem. \ref{lem-piecewise-global-phi+} because it suffices to drop the nonnegative constraints.
We omit its proof.
%\bp
%\uline{Let us first prove ``$\geq$" holds in \eqref{piecewise-global}.} Recalling that $\I(\bfz)$ denotes the support set of $\bfz$, for any $\whw = [\why; \whz]$, we have 
%\begin{align*}
%	\Xi(\whw) + \Psi_{0,\bfmu}(\whz) = & \Xi(\whw) + \sum_{i \in \I(\whz)} \mu_i \geq \min \Big\{\Xi(\bfw) + \sum_{i \in \I(\whz)} \mu_i ~\left|~ \bfz_{\OI(\whz)} = 0 \right. \Big\} \\
%	\geq & \min_{\S \subseteq [m]} \Big\{  \min \Big\{\Xi(\bfw) + \sum_{i \in \S} \mu_i ~\Big|~ \bfz_{\omS} = 0 \Big. \Big\} \Big\} 
%\end{align*}
%Minimizing both side of this inequality with respective to $\whw$ yields ``$\geq$" in \eqref{piecewise-global}.
%
%\uline{Next we prove ``$\leq$" holds in \eqref{piecewise-global}.} Given $\S \subseteq [m]$ and $\bfw = [\bfy; \bfz]$ with $\bfz_{\omS} = 0$, the following relations hold:
%\begin{align*}
%	\Xi(\bfw) + \sum_{i \in \S} \mu_i = \Xi(\bfw) + \Psi_{0,\bfmu} ( \bfz ) \geq \min_{\bfw} \Xi(\bfw) + \Psi_{0,\bfmu}(\bfz).
%\end{align*}
%Taking minimization on both side of the inequality over $\bfw = [\bfy; \bfz]$ subject to $\bfz_{\omS} = 0$, and then over $\S \subseteq [m]$, we can derive the desired conclusion.
%\ep

\begin{assumption} \label{assumption-slater} The objective function $F$ of Problem \eqref{P} is bounded below. Moreover, it holds that $\{ \bfx ~|~ \bfx \in \ri(\dom f),~ A\bfx \in \ri(\dom g),~ B\bfx = \bfb \} \neq \emptyset$. Particularly, ``\ri" in the set can be omitted if the corresponding function $f$ or $g$ is polyhedral.
\end{assumption}

%\begin{assumption} \label{assumption-slater}
%	There exists point $\bfx \in \ri(\dom f)$ such that $A\bfx \in \ri(\dom g)$ and $B\bfx = \bfb$. Particularly, ``\ri" in this statement can be omitted if the corresponding $\dom f$ or $\dom g$ is polyhedral convex.
%\end{assumption}

\begin{theorem} \label{thm-global-sol0}
Given $\Psi_\bfmu = \Psi_{0,\bfmu}$, if Assumption \ref{assumption-slater} holds, then Problem \eqref{D} has a global solution.
\end{theorem}
\bp
Assumption \ref{assumption-slater} implies that for any $\S \subseteq [r]$, the Slater condition of the following convex program holds:
\begin{align} \label{P0-S}
\inf_{\bfx} f(\bfx) + g(A\bfx) ~s.t.~ (B\bfx - \bfb)_\S = 0.
\end{align}
Its optimal value must be finite due to the lower boundedness of $F$. The dual problem of \eqref{P0-S} is as follows:
\begin{align*}
\inf_\bfw \Xi(\bfw) ~s.t.~\bfz_{\omS} = 0. 
\end{align*}
It follows from \cite[Proposition 6.4.4]{bertsekas2003convex} that the dual problem must have global solutions. Finally, by using Lem.~\ref{lem-piecewise-global-phi0}, we can conclude that \eqref{D} must have a global minimizer.
\ep

%%%%%%%%%%%%%%%%%%%%%%%%%%
\subsection{Primal-Dual Examples: Illustration} \label{Subsection-Illustration}

Our main purpose in this part is to show the proposed 
sufficient conditions
can be easily verified and may even be tight (e.g., also necessary). It implies that further improvement in general is hard.

\subsubsection{Sparse SVM problems}

We consider the SVM problem \eqref{hSVM}. 
We recall $f(\bfx) = \| \bfog\|^2/2$, $g \equiv 0$, $\Phi_\bfla(\bfu) = \Phi_{+,\bfla}(\bfu)$ and $B=- \Diag(\bfc)\overline{Q}$, $\bfb = -\bfone_r$.
The conjugate function of $f$ is given by
\[
f^*(\bfzt, \zeta_0) = \frac 12 \| \bfzt\|^2 + \delta(\zeta_0 | \{0 \}).
\]
The stationary dual problem after simplification in the form of minimization
is given by
\be \label{SVM-D} \tag{Sparse SVM}
\inf_{\bfz}    \frac 12 \| Q^\top \Diag(\bfc) \bfz \|^2 - 
\langle \bfone_r, \bfz \rangle  + \Psi_{0,\bfmu}(\bfz) \quad
\mbox{s.t.} \quad \bfc^\top \bfz =0, \ \bfz \ge 0.
\ee  
We now apply the obtained results to both the primal and dual problems.
Firstly, we note that $f$ is a convex quadratic function, and hence it is ALS. According to 
Lem.~\ref{Lemma-ALS} (i), the objective function $F$ of \eqref{hSVM} is also ALS. Moverover, for any $\bfd \in \Re^n$, we have $F_\infty(\bfd) = f_\infty(\bfd) \geq 0$ due to \eqref{asym-H} and lower boundedness of $f$. It follows from Thm.~\ref{Thm-H01} that the primal problem \eqref{hSVM} always has a global solution. 
For the dual problem, we note that $\dom g= \Re^m$ and $\dom f=  \Re^n$.
According to Thm.~\ref{thm-global-sol-D+}, if the constraints $B \bfx \leq \bfb$ are feasible, 
then there exists a global solution for the dual problem. This condition is equivalent to 
\be \label{Separablity-Condition}
\exists \  \bfx = [\bfog; \omega_0] \ \mbox{such that}\ \Diag(\bfc) \overline{Q} \bfx \ge \bfone_r .
\ee 
Interestingly, this condition is also necessary, as we prove below.

\begin{proposition} \label{Prop-SVM}
The primal SVM \eqref{hSVM} always has a global solution. 
The dual problem \eqref{SVM-D} has a global solution if and only if
condition \eqref{Separablity-Condition} holds.
\end{proposition}

\bp
We only need to prove the necessary condition part. 
Suppose the dual problem has a global solution, we prove Condition \eqref{Separablity-Condition} must hold.

We start writing the dual problem
as
\[
\inf_{\bfz}  G(\bfz) = \frac 12 \| Q^\top \Diag(\bfc) \bfz \|^2 - 
\langle \bfone_r, \bfz \rangle  + \Psi_{0,\bfmu}(\bfz) \quad
\mbox{s.t.} \quad \bfc^\top \bfz =0,~ \bfz \geq 0 .
\]
Due to the structure of $G$, it is necessary the following implication holds:
\be \label{Implication-Omega}
\bfz \in \C :=\big\{ \bfd ~ \big|~  Q^\top \Diag(\bfc) \bfd = 0, \ \bfc^\top \bfd = 0, \ \bfd \geq 0 
\big. \big \} \ \Longrightarrow \  -\langle \bfone_r, \bfz \rangle \ge 0 .
\ee
Otherwise, if there is such $\widehat{\bfz} \in \C$ such that $-\langle \bfone_r, \widehat{\bfz} \rangle < 0$. Then $\tau \widehat{\bfz} \in \C$ for $\tau >0$ (because
$\C$ is a cone) and consequently, 
\[
G( \tau \widehat{\bfz}) = -\tau \langle \bfone_r, \widehat{\bfz} \rangle + \Psi_{0,\bfmu}(\whz)  
\le -\tau \langle \bfone_r, \widehat{\bfz} \rangle + \sum_{i = 1}^r \mu_i \rightarrow - \infty \ \mbox{as} \ \tau \rightarrow + \infty
\]
contradicting that
a global solution exists. Hence, \eqref{Implication-Omega} holds.
This implication means that the linear programming below has $0$ as its optimal
objective value:
\[
\min_{\bfd} -\langle \bfone_r, \bfd \rangle, \quad \mbox{s.t.} \quad
Q^\top \Diag(\bfc) \bfd = 0, \ \bfc^\top \bfd = 0, \ \bfd \ge 0 .
\]
By the duality theory in linear programming, 
its dual problem has an optimal solution:
\be \label{LP-Dual}
\max_{\bfx = [\bfog;\omega_0]} \; 0, \quad \mbox{s.t.} \quad
\Diag(\bfc) \overline{Q} \bfx \ge \bfone_r .
\ee
Equivalently, the dual problem \eqref{LP-Dual} is feasible and condition \eqref{Separablity-Condition} holds. \hfill $\Box$

\begin{remark} \label{Remark-SVM}
The necessary and  sufficient condition \eqref{Separablity-Condition} actually is the
more familiar separability condition: there exists
$(\widehat{\bfog}, \widehat{\omega}_0)$ such that
\[
\left\{ \begin{array}{ll}
	\widehat{\bfog}^\top \bfx_i + \widehat{\omega}_0 \ge 1 & \mbox{if} \ c_i =1 \\ [0.3ex]
	\widehat{\bfog}^\top \bfx_i + \widehat{\omega}_0 \le -1 & \mbox{if} \ c_i =-1 .
\end{array} 
\right .
\]
If the data is not separable, then the dual problem has no global solutions although local solutions exist.
Our result justifies why ``regularization'' is often introduced in some existing
sparse dual SVM models. For instance, the dual model considered in \cite{zhou2021sparse} can be 
regarded as a regularized \eqref{SVM-D} with the box constraint $0 \le \bfz \le \bfga$   for some $\bfga > 0$ replacing $\bfz \ge 0$. 
It was reported therein that this box constraint improves the numerical performance of its algorithm.
Our theoretical result says that this regularization ensures existence of a global solution.
Another regularization adds a quadratic term of $\| \bfz\|^2$ to the objective.
The resulting dual problem \eqref{SVM-D} also has a global solution.
\end{remark}

%%%%%%%%%%%%%%%%%%%%%%%
\subsubsection{Energy-minimization with $\ell_0$-denoising}

We consider the problem \eqref{Energy-minimization} with the constraint $\Omega = \mbR^n_+~\mbox{or}~ \{ \bfx ~|~ 0 \leq \bfx \leq \bfga \}$, where $\gamma_i >0$ for $i \in [n]$.
We re-write the problem below as easy reference:
\be \label{Energy-Minimization-P}
\min_{\bfx} \delta(\bfx | \Omega) + \frac 12 \| A\bfx - \bfa\|^2 + \Phi_{0,\bfla}(D\bfx)
\ee 
%where
%$
% \Omega = [0, \bfu].
%$ 
%and $D_{i,i}=-1$, $D_{i,i+1} = 1$ for $i\in [n-1]$ and $D_{ij}=0$ for all other cases.
In the above case, $f(\bfx) =\delta(\bfx | \Omega)$,
$g(\bfv) = \| \bfv - \bfa\|^2/2$. When $\Omega = \{ \bfx ~|~ 0 \leq \bfx \leq \bfga \}$, we can compute $f^*(\bfzt) = \langle \max\{ \bfzt, 0 \}, \bfga \rangle$ and 
$g^*(\bfy) = \| \bfy\|^2/2 + \langle \bfa, \bfy \rangle$. 
The stationary dual problem of \eqref{Energy-Minimization-P} is given by
\be \label{Energy-Minimization-D}
\min_{\bfy, \bfz} \ \big\langle \max\{ -A^\top \bfy - D^\top \bfz, 0\}, \bfga \big\rangle  + \frac 12 \| \bfy\|^2 +  \langle \bfa, \bfy \rangle 
+ \Psi_{0,\bfmu}(\bfz).
\ee 
If $\Omega = \mbR^n_+$, then $f^*(\bfzt) = \delta(\bfzt|\mbR^n_-)$ and stationary dual problem of \eqref{Energy-Minimization-P} becomes
\be \label{Energy-Minimization-D-2}
\min_{\bfy, \bfz} \  \frac 12 \| \bfy \|^2 +  \langle \bfa, \bfy \rangle 
+ \Psi_{0,\bfmu}(\bfz), \quad \mbox{s.t.} \ \
A^\top \bfy + D^\top \bfz \ge 0 .
\ee 
Straightforward application of the obtained results leads to the following existence result.

\begin{proposition}\label{Prop-Energy}
Both the primal problem \eqref{Energy-Minimization-P} and the dual problem
\eqref{Energy-Minimization-D}/\eqref{Energy-Minimization-D-2} have global
solutions.	
\end{proposition}

\bp
For the primal problem, we note that $g$ is ALS because it is convex and quadratic, $g_\infty(A\bfd) \geq 0$ for $\bfd \in \mbR^n$ due to the lower boundedness of $g$. We see all the conditions in 
Thm.~\ref{Theorem-Primal} are met. Therefore, the primal problem has a global solution.

For the dual problem, we see $ 0 \in \ri (\dom g) = \Re^m$, $0 \in \Omega = \dom f$ and $\Omega$ is polyhedral convex, and $0 \in \Rge D$. According to Thm.~\ref{thm-global-sol0},
the dual problem \eqref{Energy-Minimization-D}/\eqref{Energy-Minimization-D-2}
has a global solution. \hfill $\Box$

For applications where the linear operator $A$ is the identity matrix ($A=I$) and $\Omega = \Re^n_+$, the dual problem \eqref{Energy-Minimization-D-2} has a simplified representation (ignoring the constant term that has no impact on optimization):
\[  %\be \label{Semismooth-Problem}
\min_{\bfz} \  \frac 12 \| ( \bfa - D^\top \bfz)_+  \|^2 
+ \Psi_{0,\bfmu}(\bfz)
\] %\ee 
The first term in the above objective function is continuously differentiable and its gradient is 
strongly semismooth.  In particular, the corresponding dual problems for the
\eqref{Edge-Denoising} and \eqref{Calcium} (both have $A=I$) are respectively given by
\[  
\min_{\bfz} \  \frac 12 \|  \bfb - D^\top \bfz  \|^2
+ \Psi_{0,\bfmu}(\bfz)  \quad \mbox{and} \quad   %\tags{Edge-D} \\
\min_{\bfz} \  \frac 12 \|  (\bfb - D^\top \bfz)_+  \|^2
+ \Psi_{+,\bfmu}(\bfz)   %\tags{Calcium-D}
\] 
%\begin{eqnarray*}  
% && \min_{\bfz} \  \frac 12 \|  \bfb - D^\top \bfz  \|^2
% + \mu \| \bfz\|_0 + \frac 12 \| \bfb\|^2 \\  %\tags{Edge-D} \\
%&& 
%\min_{\bfz \ge 0} \  \frac 12 \|  (\bfb - D^\top \bfz)_+  \|^2
%+ \mu \| \bfz\|_0 + \frac 12 \| \bfb\|^2.  %\tags{Calcium-D}
%\end{eqnarray*} 

Fast Newton-type methods can be developed for those problems. We refer to our previous papers \cite{zhang2025composite,zhang2025sparse} for the algorithmic development.
The first problem is a variant of the best subset section problem that has efficient algorithms for its global solution, see \cite{bertsimas2016best}. All illustrations above used the $\ell_2$-squared loss, which is continuously differentiable.
For the $\ell_1$ non-differentiable loss, we also get a stationary dual problem with nice structures. 

\begin{remark} \label{Remark-L1-Loss}
($\ell_1$-loss energy minimization)
Consider the primal problem:
\[
\min_{\bfx \ge 0 } \| A\bfx - \bfa\|_1 + \Psi_{0,\bfmu}(D\bfx). \label{l1-loss-energy-min}
\]
In this case, $f(\bfx) =\delta(\bfx | \mbR^n_+)$ and
$g(\bfv) = \| \bfv - \bfa\|_1$. Since $g$ is polyhedral convex and bounded below, then $g$ is ALS and $g_\infty$ is nonnegative. 
Furthermore, Thm. \ref{Theorem-Primal} ensures that \eqref{l1-loss-energy-min} has a global solution. Using the well-known fact that the conjugate function of the $\ell_1$-norm is the indicator function of $\ell_\infty$-unit ball 
($g^*(\bfy) = \langle \bfa, \bfy \rangle + \delta(\bfy \;|\; \| \cdot\|_\infty \le 1)
%\| \cdot\|_1^* = \delta(\cdot | \| \bfu\|_\infty \le 1)
$), the stationary dual problem takes the following form:
\begin{equation}\label{l1-loss-energy-min-D}
	\min_{\bfz} \ \langle \bfa, \bfy \rangle + \Psi_{0,\bfmu}(\bfz) \quad \mbox{s.t.}\ \
	A^\top \bfy + D^\top \bfz \ge 0, \ \| \bfy \|_\infty \le 1. 
\end{equation}
This is a cardinality minimization problem, which is well-surveyed in \cite{tillmann2024cardinality}. Since $\mbR^n_+$ is polyhedral convex and $\dom g = \mbR^m$, Assumption \ref{assumption-slater} holds. It follows from Thm. \ref{thm-global-sol0} that there exists a global solution of \eqref{l1-loss-energy-min-D}. 
\end{remark}

%%%%%%%%%%%%%%%%%%%%%%%%%%%%%%%%%%%
\section{Primal and Dual Optimality Analysis} \label{Section-Solution-Correspondence}

After studying global solutions of both primal and dual problems,
we now turn to establish one-to-one correspondence between solutions of them.
The analysis tools are different from the preceding  section
and require a new setup. 
%To simplify the proof,
%we only consider the case $\Phi_\bfla = \Phi_{0,\bfla}$ because the case 
%$\Phi_\bfla=\Phi_{+,\bfla}$ can be similarly treated without harming the claims we are to make in this section. 

Let us denote
\[
\left\{
\begin{array}{llll}
\bfw &:= [\bfy; \bfz] \in \mbR^{m+r}, \quad 
& \Theta(\bfx) &:= f(\bfx) + g(A\bfx)\\ [0.6ex]
Q    &:= [A; \; B] \in \mbR^{(m+r)\times n}, \quad
& \Xi(\bfw) &:= f^*(-Q^\top \bfw) + g^*(\bfy) + \langle \bfb, \bfy \rangle.
\end{array} 
\right .
\]

%\begin{align*}
%	\bfw := [\bfy; \bfz] \in \mbR^{m+r} ~~&\mbox{and}~~ Q := [A; \; B] \in \mbR^{(m+r)\times n}, \\
%	\Psi_\bfmu(\bfx) := f(\bfx) + g(A\bfx) ~~&\mbox{and}~~ \Phi_\bfla(\bfw) := f^*(-Q^\top \bfw) + g^*(\bfy),
%\end{align*} 
Problems \eqref{P} and \eqref{D} are restated as follows:
\begin{align*}
\mbox{(P)} \quad 	\min_{\bfx }\; F(\bfx) =\Theta(\bfx) + \Phi_\bfla(B\bfx - \bfb) \quad \mbox{and} \quad 
\min_{\bfw}\; G(\bfw) = \Xi(\bfw) + \Psi_\bfmu(\bfz), \qquad \mbox{(D)}
\end{align*}
where $\Phi_\bfla \in \{ \Phi_{0,\bfla}, \Phi_{+,\bfla} \}$ and $\Psi_\bfmu \in \{ \Psi_{0,\bfmu}, \Psi_{+,\bfmu} \}$. We note that the functions $\Theta(\bfx)$ and $\Xi(\bfw)$ are respectively the
convex part of the primal and the dual problems.
In our solution correspondence characterization below, we will match the solutions that admit equal function
values for the convex parts, i.e., $\Theta(\bfx^*) = - \Xi(\bfw^*)$.

We start with characterizing the optimality of \eqref{P} and \eqref{D} by their stationary points. The definitions are based on Fermat's rule \cite[Thm.~10.1]{RockWets98} formulated via the limiting subdifferential.
Similar definitions have also been given in \cite{bolte2018nonconvex,cai2022developments,themelis2018forward}.

\begin{definition} \label{def-stationary-P}
We say $\bfx^*$ is a stationary point of \eqref{P} if it satisfies 
\begin{equation} \label{stationary-P}
	0 \in \partial f(\bfx^*) + A^\top \partial g(A\bfx^*) + B^\top \partial \Phi_\bfla(B \bfx^* - \bfb).
\end{equation}
Moreover, $\bfw^*:= [\bfy^*;\bfz^*]$ is a stationary point of \eqref{D} if it satisfies
\begin{align} \label{stationary-D}
	0 \in - \left( \begin{array}{c}
		A \\
		B
	\end{array} \right) \partial f^*( -Q^\top \bfw^* ) + \left( \begin{array}{c}
		\partial g^*(\bfy^*) \\
		\partial \Psi_\bfmu(\bfz^*)
	\end{array} \right). 
\end{align}
\end{definition}

%\begin{definition}
%	We say $\bfw^*:= [\bfy^*;\bfz^*]$ is a stationary point of \eqref{D} if it satisfies
%	\begin{align} \label{stationary-D}
%		0 \in - \left( \begin{array}{c}
	%			A \\
	%			B
	%		\end{array} \right) \partial f^*( -Q^\top \bfw^* ) + \left( \begin{array}{c}
	%			\partial g^*(\bfy^*) \\
	%			\partial \mu \card{\bfz^*}
	%		\end{array} \right). 
%	\end{align}
%\end{definition}

The road of our proofs goes like this. 
For each problem of (P) and (D), we prove that a local solution
is a stationary point. We further prove that a stationary point
of one problem corresponds to a stationary point of the other. 
Through this, we prove the one-to-one correspondence of the local solutions of (P) and (D). 
For such correspondence to hold, we will need a weak form of 
the Slater condition.

%\subsection{Optimality Condition of Primal Problem}
\subsection{Equivalence of local minimizers and stationary points of \eqref{P}} 

\subsubsection{\texorpdfstring{Problem~\eqref{P} with $\Phi_{\bfla}=\Phi_{+,\bfla}$}
{Problem (P) with Phi_lambda = Phi_+,lambda}} \label{sec-P-phi+}
%\subsubsection{The case of \eqref{P} with } \label{sec-P-phi+}
In classic variational analysis theory, establishing \eqref{stationary-P} as a necessary optimality condition of \eqref{P} often requires some regularity conditions (see \cite[Theorems.~10.1, 10.6, and Corollary 10.9]{RockWets98}). However, we will see that our analysis can be simplified by leveraging the convexity of $f$ and $g$, and piecewise-constant structure of the cardinality function $\Phi_{+,\bfla}$. Given a reference point $\bfx^* \in \mbR^n$, we denote $\J_* := \I(B\bfx^* - \bfb)$ (the support set of $B\bfx^* - \bfb$, hence $\OJ_*$ is the set of indices of $i$ with $(B \bfx^* - \bfb)_i=0$) and define the following convex program associated with $\J_*$:
\begin{align} \tag{${\rm P}_+$} \label{P+}
\min_{\bfx} \Theta(\bfx) \ \ s.t. \ \ (B\bfx - \bfb)_{\OJ_*} \leq 0.
\end{align}
We say $\bfx$ is a KKT (Karush-Kuhn-Tucker) point of \eqref{P+}, then there exists a multiplier $\bfz_{\OJ_*} \in \mbR^{|\OJ_*|}_+$ such that $(\bfx,\bfz_{\OJ_*})$ satisfies the following KKT system
\begin{align} \label{KKT-P+}
\left\{	\begin{aligned}
	& 0 \in \partial f(\bfx) + A^\top \partial g(A\bfx) + B_{\OJ_*:}^\top \bfz_{\OJ_*}, \\
	& 0 \leq \bfz_{\OJ_*} \perp (B\bfx - \bfb)_{\OJ_*} \leq 0.
\end{aligned} \right.
\end{align} 
As \eqref{P+} is a convex problem, a KKT point must be a global minimizer, and the converse conclusion is guaranteed under the generalized Slater condition
\cite{bertsekas2003convex,boyd2004convex} of \eqref{P+}:
\begin{align} \label{slater-P+}
\{ \bfx \ | \ \bfx \in \ri(\dom f),~ A\bfx \in \ri(\dom g),~ (B\bfx - \bfb)_{\OJ_*} \leq 0 \} \neq \emptyset,
\end{align} 
In particular, if function $f$ (resp. $g$) is polyhedral,
the associated ``{\ri}" in \eqref{slater-P+} can be dropped and the condition is merely the feasibility condition.
\begin{lemma} \label{lem-relation-P-P+}
Given $\bfx^* \in \mbR^n$ and $\Phi_\bfla = \Phi_{+,\bfla}$, the following assertions hold: 
%	The following assertions hold for problems \eqref{P} and \eqref{P*}:

(i) $\bfx^*$ is a stationary point of \eqref{P} if and only if it is a KKT point of \eqref{P+}.

(ii) $\bfx^*$ is a local minimizer of \eqref{P} if and only if it is a global minimizer of \eqref{P+}.
\end{lemma}
\bp
(i) This is due to the structure of the limiting subdifferential $\partial \Phi_{+,\bfla}$ in \eqref{subdiff-phi+} and the definition of $\J_*$.
We can verify that a point $\bfx^*$ satisfies \eqref{stationary-P} if and only if there exists vector $\bfz^*$ with $\bfz^*_{\J_*} = 0$ and $\bfz^*_{\OJ_*} \in \mbR^{|\OJ_*|}$ such that \eqref{KKT-P+} holds for $(\bfx^*,\bfz^*_{\OJ_*})$.

(ii) Let us denote the feasible region of \eqref{P} and \eqref{P+} 
by $\F$ and $\F_*$ respectively.
%as $\F : = \{ \bfx \ | \ \bfx \in \dom f,~ A\bfx \in \dom g \}$ and $\F_* = \{ \bfx \ | \ \bfx \in \dom f,~ A\bfx \in \dom g, ~ (B\bfx)_{\OG_*} = 0 \}$ respectively.

``$\Rightarrow$" If $\bfx^*$ is a local minimizer of \eqref{P}, then there exists radius $\epsilon_1 >0$ such that 
\begin{align} \label{loc-min-P+}
\Theta(\bfx) + \Phi_{+,\bfla}(B\bfx - \bfb) \geq \Theta(\bfx^*) + \Phi_{+,\bfla}(B\bfx^* - \bfb),~ \forall \ \bfx \in \F \cap \N(\bfx^*, \epsilon_1).
\end{align}
When $\epsilon_1$ is small enough, we also have 
\begin{align} \label{sgn-J}
\sgn(B\bfx - \bfb)[i] = \sgn(B\bfx^* - \bfb)[i]~\mbox{for}~ i \in \J_*
\end{align}
Moreover, for $\bfx \in \F_* \cap \N(\bfx^*, \epsilon_1)$, we have $(B\bfx - \bfb)_{\OJ_*} \leq 0$. This together with \eqref{sgn-J} imply $\Phi_{+,\bfla}(B\bfx - \bfb) = \Phi_{+,\bfla}(B\bfx^* - \bfb)$. Combining this with \eqref{loc-min-P+}, we can derive
\begin{align*}
\Theta(\bfx) \geq \Theta(\bfx^*),~ 
\forall \ \bfx \in \F_* \cap \N(\bfx^*, \epsilon_1).
\end{align*}
That is, $\bfx^*$ is a local minimizer of \eqref{P+}. Since \eqref{P+} is a convex program, $\bfx^*$ must be its global minimizer. 

``$\Leftarrow$" If $\bfx^*$ is a global minimizer of \eqref{P+}, then we have
%$\Psi_\bfmu(\bfx) \geq \Psi_\bfmu(\bfx^*)$ for all $\bfb \in \F_*$.
\begin{align} \label{global-min-P*}
\Theta(\bfx) \geq \Theta(\bfx^*), ~ \forall \ \bfx \in \F_*.
\end{align}
Let us take a sufficiently small radius $\epsilon_2$ such that \eqref{sgn-J} and the following formula hold for any $\bfx \in \F \cap \N(\bfx^*, \epsilon_2)$:
\begin{align}
%	&\I(B\bfx - \bfb) \supseteq \Ga_* := \J_*, \label{T-supseteq-T*} \\
\Theta(\bfx) \geq \Theta(\bfx^*) - \min_{i \in [r]}\lambda_i/2, \label{lsc-property}
\end{align}
where \eqref{lsc-property} 
follows from the lower semicontinuity of $\Theta(\cdot)$. Now let us consider the following two cases.

Case I: Taking $\bfx \in \F_* \cap \N(\bfx^*,\epsilon_2)$, then \eqref{sgn-J} implies $\Phi_{+,\bfla}(B\bfx - \bfb) \geq \Phi_{+,\bfla}(B\bfx^* - \bfb)$. 
%Together with the fact $\Psi_\bfmu(\bfx) \geq \Psi_\bfmu(\bfx^*)$ for all $\bfb \in \F_*$
Combining this with \eqref{global-min-P*} yields
\begin{align*}
\Theta(\bfx) + \Phi_{+,\bfla}(B\bfx - \bfb) \geq \Theta(\bfx^*) + \Phi_{+,\bfla}(B\bfx^* - \bfb),~ \forall \ \bfx \in \F_* \cap \N(\bfx^*, \epsilon_2).
\end{align*} 

Case II: Taking $\bfx \in (\F \backslash \F_*) \cap \N(\bfx^*, \epsilon_2)$, then $(B\bfx^* -\bfb)_{\OJ_*} > 0$. This together with \eqref{sgn-J} indicate $\Phi_{+,\bfla}(B\bfx - \bfb) \geq \Phi_{+,\bfla}(B\bfx^* - \bfb) + \min_{i \in [r]} \lambda_i$. Using \eqref{lsc-property}, we can obtain the following inequality for any $\bfx \in (\F\backslash\F_*) \cap \N(\bfx^*, \epsilon_2)$:
\begin{align*}
\Theta(\bfx) + \Phi_{+,\bfla}(B\bfx - \bfb) \geq \Theta(\bfx^*) + \Phi_{+,\bfla}(B\bfx^* - \bfb) + \min_{i \in [r]} \lambda_i/2.
\end{align*}
Overall, the above two cases imply $\bfx^*$ is a local minimizer of \eqref{P}.
\ep

This further helps to establish the optimality condition of \eqref{P} with $\Phi_\bfla = \Phi_{+,\bfla}$.
\begin{theorem} \label{thm-optimality-P+}
About problem \eqref{P} with $\Phi_\bfla = \Phi_{+,\bfla}$, we have:

(i) A stationary point is a local minimizer.

(ii) If $\bfx^*$ is a local minimizer and Slater condition \eqref{slater-P+} holds, then it is a stationary point.	
\end{theorem}
\bp
As we have stated just before Lem.~\ref{lem-relation-P-P+}, a KKT point of \eqref{P+} is a global minimizer of \eqref{P+} and the converse is true when \eqref{slater-P+} holds. Then by using Lem.~\ref{lem-relation-P-P+}, we can derive the desired conclusion.
\ep

It is worth mentioning that the recent works \cite{cui2023minimization,hananalysis} use pseudo B-stationary point and epi-stationary point to characterize the local solutions of problems with cardinality functions. In our paper, the stationary point based on limiting subdifferential is adopted because we need to use the stationary duality \eqref{stationary-duality} and \eqref{stationary-duality+} to establish the correspondence of local solutions to \eqref{P} and \eqref{D}.

%{compare with cui, we use stationary point based on subdifferential, this helps to establish primal-dual correspondence with the aid of subdifferential conjugacy}

\subsubsection{\texorpdfstring{Problem~\eqref{P} with $\Phi_{\bfla}=\Phi_{0,\bfla}$}
{Problem (P) with Phi_lambda = Phi_0,lambda}} \label{P-Phi=Phi0}
%\subsubsection{The case of \eqref{P} with } \label{P-Phi=Phi0}
We can notice that $\Phi_{0,\bfla}(\bfu) = \Phi_{+,\bfla}(\bfu) + \Phi_{+,\bfla}(-\bfu)$. Therefore, denoting $\wtB:=[B;-B]$ and $\wtb:= [\bfb;-\bfb]$, Problem \eqref{P} with $\Phi_\bfla = \Phi_{0,\bfla}$ can be equivalently reformulated as 
\begin{align} \label{P-phi0-eq}
\min_\bfx \Theta(\bfx) + \Phi_{+,\bfla}\big(\wtB \bfx - \wtb\big)
\end{align}
In this way, we can utilize the results in Section \ref{sec-P-phi+} to establish the optimality condition of Problem \eqref{P} with $\Phi_\bfla = \Phi_{0,\bfla}$.
%Given a reference point $\bfx^* \in \mbR^n$, let us denote $\Ga^*:= \I(B)$

With the special structure of $\wtB$ and $\wtb$ in \eqref{P-phi0-eq}, given a reference point $\bfx^* \in \mbR^n$ and index set $\J_* = \I(B\bfx^* - \bfb)$, the Slater condition \eqref{slater-P+} can be simplified as follows:
\begin{align} \label{slater-P0}
\{ \bfx \ | \ \bfx \in \ri(\dom f),~ A\bfx \in \ri(\dom g),~ (B\bfx - \bfb)_{\OJ_*} = 0 \} \neq \emptyset,
\end{align} 
where ``\ri" in the above set can be omitted if the associated function $f$ or $g$ is polyhedral convex. Next, we characterize the optimality condition of Problem \eqref{P} with $\Phi_\bfla = \Phi_{0,\bfla}$ by its stationary point.
\begin{theorem} \label{thm-optimality-P}
About problem \eqref{P} with $\Phi_\bfla = \Phi_{0,\bfla}$, we have:

(i) A stationary point is a local minimizer.

(ii) If $\bfx^*$ is a local minimizer and Slater condition \eqref{slater-P+} holds, then it is a stationary point.	
\end{theorem}
\bp
According to Thm. \ref{thm-optimality-P+}, the assertions (i) and (ii) hold for Problem \eqref{P-phi0-eq}. It suffices to show that the stationary points of \eqref{P-phi0-eq} and \eqref{P} with $\Phi_\bfla = \Phi_{0,\bfla}$. Let $\bfx^*$ be a stationary point of Problem \eqref{P-phi0-eq}. It satisfies
\begin{equation} \label{stat-phi0-eq} 
0 \in \partial f(\bfx^*) + A^\top \partial g(A\bfx^*) + \wtB^\top \partial \Phi_{+,\bfla}\big(\wtB \bfx^* - \wtb \big).
\end{equation}
By the representation of $\partial \Phi_{0,\bfla}$ and $\partial \Phi_{+,\bfla}$, we have $\partial \Phi_{0,\bfla}(\bfu) = \partial \Phi_{+,\bfla}(\bfu) - \partial \Phi_{+,\bfla}(-\bfu)$ for any $\bfu \in \mbR^r$. This together with \eqref{stat-phi0-eq} lead to the conclusion.
\ep

\subsection{Equivalence of local minimizers and stationary points of \eqref{D}}

%\subsubsection{The case of \eqref{D} with $\Phi_\bfla = \Phi_{+,\bfla}$}
\subsubsection{\texorpdfstring{The case of \eqref{D} with $\Psi_\bfmu = \Psi_{+,\bfmu}$}{The case of D with bfla = Psi_+,bfla}}
This is a parallel development for the dual problem. To establish the equivalence of local minimizers and stationary points of \eqref{D}, we also need to introduce a convex problem associated with certain index set. Given a reference point $\bfw^* = [\bfy^*;\bfz^*]$, let us denote $T_*:= \I(\bfz^*)$ and consider the following convex optimization:
\begin{align}  \tag{${\rm D}_+$} \label{D+}
\min_{\bfw = [\bfy; \bfz] }\; 
\Xi(\bfw) \quad \mbox{s.t.} \ \bfz_{\oT_*} = 0,~ \bfz_{T_*} \geq 0.
\end{align}
If $\bfw$ is a KKT point of \eqref{D+}, there exists $\bfu \in \mbR^r$ with $\bfu_{T_*} \in \mbR^{|T^*|}_-$ and $\bfu_{\oT_*} \in \mbR^{|\oT_*|}$ such that 
\begin{align} \label{KKT-D+}
\left\{ \begin{aligned}
	& 0 \in - \left( \begin{array}{c}
		A \\
		B
	\end{array} \right) \partial f^*( -Q^\top \bfw ) + \left( \begin{array}{c}
		\partial g^*(\bfy) \\
		\bfu
	\end{array} \right), \\
	& 0 \leq \bfz_{T_*} \perp \bfu_{T_*} \leq 0.
\end{aligned} \right.
\end{align}
Since \eqref{D+} is convex, a KKT point must be a global minimizer, and the converse conclusion holds under the following Slater condition:
\begin{align} \label{slater-D+}
\{ \bfw = [\bfy;\bfz] | -Q^\top\bfw \in \ri(\dom f^*), ~ \bfy \in \ri(\dom g^*), ~\bfz_{\oT_*} = 0,~ \bfz_{T_*} \geq 0 \} \neq \emptyset,
\end{align}
where ``\ri" in this condition can be dropped when the associated function $f^*$ or $g^*$ is polyhedral convex. Problems \eqref{D} and \eqref{D+} has close relation, which is summarized in the following lemma.
\begin{lemma} \label{lem-relation-D-D+}
Given $\bfw^* = [\bfy^*;\bfz^*]$ and $\Psi_\bfmu = \Psi_{+,\bfmu}$, the following assertions hold: 
%	The following assertions hold for problems \eqref{P} and \eqref{P*}:

(i) $\bfw^*$ is a stationary point of \eqref{D} if and only if it is a KKT point of \eqref{D+}.

(ii) $\bfw^*$ is a local minimizer of \eqref{D} if and only if it is a global minimizer of \eqref{D+}.
\end{lemma}
\bp
(i) By the representation of $\partial \Psi_{+,\bfmu}$ and the 
definition of $T_*$, we have 
\begin{align*}
0 \leq \bfz^*_{T_*} \perp \bfu^*_{T_*} \leq 0 \iff \bfu^* \in \partial \Psi_{+,\bfmu}(\bfz^*). 
\end{align*}
Then comparing \eqref{KKT-D+} and \eqref{stationary-P}, we can arrive at the desired conclusion.

$\bfu_{T_*} = 0$ and $\bfu_{\oT_*} \in \mbR^{|\oT_*|}$. Therefore, we have $\bfu \in  \partial\Psi_{+,\bfmu}(\bfz) $ in \eqref{KKT-D+}, which means (i) is true.

(ii) Let us denote the feasible regions of \eqref{D} and \eqref{D+} as $\G$ and $\G_*$ respectively.
%Let us denote the feasible regions of \eqref{D} and \eqref{D+} as $\G:=\{ \bfw = [\bfy;\bfz] \ | \ -Q^\top \bfw \in \dom f^*,~ \bfy \in \dom g^* \}$ and $\G_* := \{ \bfw = [\bfy;\bfz] \ | \ -Q^\top \bfw \in \dom f^*,~ \bfy \in \dom g^*,~ \bfz_{\oT_*} = 0 \}$.

``$\Rightarrow$" If $\bfw^*$ is a local minimizer of \eqref{D}, then there exists $\epsilon_3 > 0$ such that 
\begin{align*}
\Xi(\bfw) + \Psi_{+,\bfmu}(\bfz) \geq \Xi(\bfw^*) + \Psi_{+,\bfmu}(\bfz^*), \ \forall \ \bfw \in \G \cap \N(\bfw^*,\epsilon_3).
\end{align*} 
If we further take $\bfw \in \G_* \cap \N(\bfw^*,\epsilon_3)$, then $\bfz_{\oT_*} = 0$ and therefore $\Psi_{+,\bfmu}(\bfz) \leq \Psi_{+,\bfmu}(\bfz^*)$ by the definition of $T_*$. Combining this with the above inequality yields
\begin{align*}
\Xi(\bfw) \geq \Xi(\bfw^*), \ \forall \ \bfw \in \G_* \cap \N(\bfw^*, \epsilon_3).
\end{align*}
Considering that \eqref{D+} is a convex program, $\bfw^*$ must be a global minimizer of \eqref{D+}.

``$\Leftarrow$" If $\bfw^*$ is a global minimizer of \eqref{D+}, then we have
\begin{align} \label{global-min-D0}
\Xi(\bfw) \geq \Xi(\bfw^*), \ \forall \ \bfw \in \G_*.
\end{align} 
Now let us take a sufficiently small radius $\epsilon_4 >0$ such that the following formulas hold for any $\bfw \in \G \cap \N(\bfw^*, \epsilon_4)$:
\begin{align}
& \I(\bfz) \supseteq T_* := \I(\bfz^*), \label{Iz-supset-Iz*} \\
& \Xi(\bfw) \geq \Xi(\bfw^*) - \min_{i \in [r]}\mu_i/2, \label{lsc-dual}
\end{align}
where the second line above follows from the lower semicontinuity of $\Xi$. Next we consider the following two cases.

Case I: If $\bfw \in \G_* \cap \N(\bfw^*, \epsilon_4)$, then \eqref{Iz-supset-Iz*} implies $\Psi_{+,\bfmu}(\bfz) \geq \Psi_{+,\bfmu}(\bfz^*)$. Therefore, it follows from \eqref{global-min-D0} that
\begin{align*}
\Xi(\bfw) + \Psi_{+,\bfmu}(\bfz) \geq \Xi(\bfw^*) + \Psi_{+,\bfmu}(\bfz^*), \ \forall \ \bfw \in \G_* \cap \N(\bfw^*, \epsilon_4).
\end{align*} 

Case II: If $\bfw \in (\G \backslash \G_*) \cap \N(\bfw^*, \epsilon_4)$, then $\bfz_{\oT_*} \neq 0$, and therefore \eqref{Iz-supset-Iz*} implies $\Psi_{+,\bfmu}(\bfz) \geq \Psi_{+,\bfmu}(\bfz^*) + \min_{i \in [r]}\mu_i$. Combining this with \eqref{lsc-dual}, we can obtain
\begin{align*}
\Xi(\bfw) + \Psi_{+,\bfmu}(\bfz) \geq \Xi(\bfw^*) + \Psi_{+,\bfmu}(\bfz^*)+ \min_{i \in [r]}\mu_i/2, \ \forall \ \bfw \in (\G \backslash \G_*) \cap \N(\bfw^*, \epsilon_4).
\end{align*} 
Finally, the above two cases indicate $\bfw^*$ is a local minimizer of \eqref{D}.\ep
\begin{theorem} \label{thm-optimality-D+}
About Problem \eqref{D} with $\Psi_\bfmu = \Psi_{+,\bfmu}$, we have:

(i) A stationary point is a local minimizer.

(ii) If $\bfw^*$ is a local minimizer and Slater condition \eqref{slater-D+} holds, then it is a stationary point.	
\end{theorem}
\bp
As we have analyzed before Lem.~\ref{lem-relation-D-D+}, a KKT point of \eqref{D+} is a global minimizer of \eqref{D+} and the converse is true when Slater condition \eqref{slater-D+} holds. Then by using Lem.~\ref{lem-relation-D-D+}, we can obtain the desired conclusion.
\ep

\subsubsection{\texorpdfstring{The case of \eqref{D} with $\Psi_\bfmu = \Psi_{0,\bfmu}$}{The case of D with Psi_bfla = Psi_0,bfla$} } \label{sec-D-phi+}
%\subsubsection{The case of \eqref{D} with } \label{sec-D-phi+}
To establish the equivalence of a local minimizer and stationary point of \eqref{P} in this case, we just need to follow a similar procedure to that in Section \ref{P-Phi=Phi0}. The key distinction from the previous section is the introduction of a new convex program:
\begin{align}  \tag{${\rm D}_0$} \label{D0}
\min_{\bfw = [\bfy; \bfz] }\; 
\Xi(\bfw) \quad \mbox{s.t.} \ \bfz_{\oT_*} = 0.
\end{align}
A point $\bfw$ is a KKT point of \eqref{D0} if there exists $\bfu \in \mbR^r$ with $\bfu_{T_*} = 0$ such that 
\begin{align} \label{KKT-D0}
0 \in - \left( \begin{array}{c}
	A \\
	B
\end{array} \right) \partial f^*( -Q^\top \bfw ) + \left( \begin{array}{c}
	\partial g^*(\bfy) \\
	\bfu
\end{array} \right). 
\end{align}
Since \eqref{D0} is convex, its KKT point must be its global minimizer and the converse is true when the following generalized Slater condition holds:
% (see \cite[Theorems 23.5, 23.8, and 23.9]{rockafellar1970convex}):
\begin{align} \label{slater-D0}
\{ \bfw = [\bfy;\bfz] \ | \ -Q^\top\bfw \in \ri(\dom f^*), ~ \bfy \in \ri(\dom g^*), ~\bfz_{\oT_*} = 0 \} \neq \emptyset,
\end{align}
where the ``\ri" can be omitted if the corresponding function $f^*$ or $g^*$ is polyhedral. For \eqref{Energy-minimization} with $\Omega = \{ \bfx ~|~ 0 \leq \bfx \leq \bfga\}$ or $\mbR^n_+$, the stationary dual problems are \eqref{Energy-Minimization-D} and \eqref{Energy-Minimization-D-2} respectively. In these two cases, we can compute $\dom f^* = \mbR^n$ or $\mbR^n_-$, and $\dom g^* = \mbR^m$. Therefore, these two stationary dual problems satisfy \eqref{slater-D0}. 
%We need to stress that the Slater condition can be easily fulfilled. 
%For example, if $g$ is strongly convex (e.g., $g(\bfu) = \| \bfu\|^2$), then $\dom g^* = \Re^m$. If $f$ is polyhedral, so is $f^*$ \cite{rockafellar1970convex}. 
%In this situation covering the problem \eqref{Energy-minimization} with $\Omega$ being polyhedral), the Slater conditions reduces to the feasibility of \eqref{D0}.

%Next let us demonstrate the close relation between \eqref{D} and \eqref{D0}.

Next, we give the relationship between Problems \eqref{D} and \eqref{D0}.

\begin{lemma} \label{lem-relation-D-D0}
Given $\bfw^* = [\bfy^*;\bfz^*]$ and $\Psi_\bfmu = \Psi_{0,\bfmu}$, the following assertions hold: 
%	The following assertions hold for problems \eqref{P} and \eqref{P*}:

(i) $\bfw^*$ is a stationary point of \eqref{D} if and only if it is a KKT point of \eqref{D0}.

(ii) $\bfw^*$ is a local minimizer of \eqref{D} if and only if it is a global minimizer of \eqref{D0}.
\end{lemma}
\bp
By the definition of $T_*$ and representation of $\partial \Psi_{0,\bfmu}$, we have
\begin{align*}
\bfu^*_{T_*} = 0 \iff \bfu^* \in \partial \Psi_{0,\bfmu}(\bfz^*).
\end{align*}
Comparing \eqref{KKT-D0} and \eqref{stationary-D}, we can prove (i). To prove assertion (ii), we follow the same argument as in Lem.~\ref{lem-relation-D-D+}, with the nonnegativity constraints removed. The details of the proof are omitted for brevity.
\ep

Finally, by using Lem. \ref{lem-relation-D-D0} and Slater condition \eqref{slater-D0}, we can characterize the local minimizer of Problem \eqref{D} with $\Psi_\bfmu = \Psi_{+,\bfmu}$ by its stationary point.
\begin{theorem} \label{thm-optimality-D}
About problem \eqref{D} with $\Psi_\bfmu = \Psi_{0,\bfmu}$, we have:

(i) A stationary point is a local minimizer.

(ii) If $\bfw^*$ is a local minimizer and Slater condition \eqref{slater-D0} holds, then it is a stationary point.	
\end{theorem}

\subsection{Correspondence of primal and dual solutions}

\begin{theorem} \label{thm-correspondence}
About Problems \eqref{P} with $\Phi_\bfla = \Phi_{+,\bfla}$ (resp. $\Phi_\bfla = \Phi_{0,\bfla}$) and \eqref{D} with $\Psi_\bfmu = \Psi_{+,\bfmu}$ (resp. $\Psi_\bfmu = \Psi_{0,\bfmu}$), we have

(i) If $\bfx^*$ is a local minimizer of \eqref{P} and Slater condition \eqref{slater-P+} (resp. \eqref{slater-P0}) holds, then there exist $\bfy^* \in \partial g(A\bfx^*)$ and $\bfz^* \in \partial \Phi_\bfla(B\bfx^* - \bfb)$ such that $\bfw^* = [\bfy^*; \bfz^*]$ is a local minimizer of \eqref{D}. Moreover, it holds that $\Theta(\bfx^*) = -\Xi(\bfw^*)$.

(ii) If $\bfw^*$ is a local minimizer of \eqref{D} and Slater condition \eqref{slater-D+} (resp. \eqref{slater-D0}), then there exists $\bfx^* \in \partial f^*(-Q^\top \bfw^*)$ such that $\bfx^*$ is a local minimizer of \eqref{P}. Moreover, it holds that $\Theta(\bfx^*) = -\Xi(\bfw^*)$.
\end{theorem}
\bp We only consider the case of \eqref{P} with $\Phi_\bfla = \Phi_{+,\bfla}$ and \eqref{D} with $\Psi_\bfmu = \Psi_{+,\bfmu}$ because the other case can be proved following a similar procedure.

(i) By Thm. \ref{thm-optimality-P+} (ii), $\bfx^*$ satisfies \eqref{stationary-P}, which means that there exist $\bfy^* \in \partial g(A\bfx^*)$ and $\bfz^* \in \partial \Phi_{+,\bfla}(B\bfx^* -\bfb)$ such that $-Q^\top \bfw^* \in \partial f(\bfx^*)$. Then by property \eqref{stationary-duality+}, we can obtain
\begin{align} \label{stationary-point-psi}
A\bfx^* \in \partial g^*(\bfy^*),~ B\bfx^*\in \partial \Psi_{+,\bfmu}(\bfz^*),~ \mbox{and}~ \bfx^* \in \partial f^*(-Q^\top \bfw^*).
\end{align}
Therefore, $\bfw^*$ satisfies \eqref{slater-D0} and it is a local minimizer of \eqref{D} with $\Psi_\bfmu = \Psi_{+,\bfmu}$ by Thm. \ref{thm-optimality-D+} (i). Moreover, it follows from \eqref{stationary-point-psi}, \eqref{conjugate-subgradient}, and \eqref{subdiff-psi+} that 
\begin{align} \label{local-duality}
\Theta(\bfx^*) + \Xi(\bfw^*) = \langle - Q^\top \bfw^*, \bfx^* \rangle + \langle \bfy^*, A\bfx^* \rangle = - \langle B\bfx^*, \bfz^*  \rangle = 0 
\end{align}
where the last equality holds because $B\bfx^*\in \partial \Psi_{+,\bfmu}(\bfz^*)$ and \eqref{subdiff-psi+} implies $\langle B\bfx^*, \bfz^*  \rangle = 0$.

(ii) It follows from Thm. \ref{thm-optimality-D+} (ii) that $\bfw^*$ satisfies \eqref{stationary-D}. Then there exists $\bfx^* \in \partial f^*(-Q^\top \bfw^*)$ such that $A\bfx^* \in \partial g^*(\bfy^*)$ and $B \bfx^* \in \partial \Psi_{+,\bfmu}(\bfz^*)$. Then by property \eqref{stationary-duality+}, we have
\begin{align} \label{stationary-phi}
\bfy^* \in \partial g(A\bfx^*),~ \bfz^* \in \partial \Phi_{+,\bfla}(B\bfx^* -\bfb), ~\mbox{and}~ -Q^\top \bfw^* \in \partial f(\bfx^*).
\end{align}
This means that $\bfx^*$ satisfies \eqref{stationary-P} and it is also a local minimizer of \eqref{P} by Thm. \ref{thm-optimality-P+} (i). Finally, through the use of \eqref{conjugate-subgradient}, \eqref{subdiff-phi+}, we can derive \eqref{local-duality}. \ep

\begin{remark}
From the proof of Thm.~\ref{thm-correspondence}, we can also derive the correspondence of stationary points of \eqref{P} and \eqref{D} without assuming Slater conditions. 
Thm. \ref{thm-correspondence} provides a theoretical guarantee for the validity of \eqref{D}. 
In particular, a solution of \eqref{P} can be exactly recovered by solving \eqref{D} in certain cases. For example, when $f$ is strongly convex (hence its conjugate $f^*$ is continuously differentiable) and a stationary point $\bfw^*$ of \eqref{D} is given, we can calculate a stationary point of \eqref{P} by $\bfx^* = \nabla f^*(-Q^\top \bfw^*)$. 
When $f$ is merely convex and lsc, the computation of $\bfx^*$ can be
constructed through algorithms. 
We will explore this direction of research in future.

%	However, this becomes challenging when $f$ is merely convex and lsc. In the next section, we will design a proximal point algorithms for solving \eqref{D}. Meanwhile, it generates a primal sequence whose accumulation points are stationary points of \eqref{P} with mild assumptions. 
\end{remark}

It is also important to note the limitation of Thm.~\ref{thm-correspondence}, which only states 
corresponding relationship between local solutions. 
In other words, if $\bfx^*$ is a global solution of the primal problem \eqref{P}, Thm.~\ref{thm-correspondence} says that its corresponding dual solution $\bfw^*$ is only a local 
solution. Fortunately, the local solution set of the dual problem does not change with the choice of  
%Next, let us discuss the relationship between the local and global solutions of \eqref{D} under the Assumption \ref{assumption-slater}. An intriguing observation from Thm. \ref{thm-optimality-D} is that the local solution set of \eqref{D} is independent of 
regularization parameter $\bfmu$ (i.e., independent of $\bfmu$).
This is because,  for $\Psi_\bfmu \in \{ \Psi_{0,\bfmu}, \Psi_{+,\bfmu} \}$,
the set $\partial\Psi_\bfmu (\bfz)$ is identical for all $\bfmu > 0$. However, the set of global solution of \eqref{D} changes with $\bfmu$. 
Given a local minimizer $\bfw^*$ of \eqref{D}, we will show that it becomes a global minimizer of \eqref{D} for appropriately selected $\bfmu$. 

Let us first consider the case of $\Psi_\bfmu = \Psi_{+,\bfmu}$. Lem. \ref{lem-relation-D-D+} indicates that a local minimizer $\bfw^*$ belongs to the following set of solution with 
$T_*= \I(\bfz^*)$:
\begin{align} \label{omega+} 
\Omega_+:= \argmin_{\bfw = [\bfy; \bfz] }\; 
\Xi(\bfw) \quad \mbox{s.t.} \ \bfz_{\oT_*} = 0,~\bfz_{T_*} \geq 0.
\end{align}
We assume that $\bfw^*$ is the point in $\Omega_+$ with the smallest cardinality on $\bfz$ variable:  
\begin{align} \label{card-min}
\bfw^* \in \argmin_\bfw \{ \| \bfz \|_0 ~|~ \bfw \in \Omega_+ \}.
\end{align}
Recalling Lem. \ref{lem-piecewise-global-phi+}, we reformulate Problem \eqref{D} as minimization of finitely many convex programs with respect to $\S \subseteq [r]$:
\begin{align} \label{convex-S}
\min_\bfw \Xi(\bfw) ~s.t.~ \bfz_{\omS} = 0,~ \bfz_\S \geq 0. \big. 
\end{align}
To proceed, we consider three cases of the above convex programs by taking $\S = T_*$, $\S \subsetneq T_*$, and $\S \nsubseteq T_*$, which yield the following three optimal values
\begin{align} \label{para-eta}
\begin{aligned}  
	& \eta_0 :=   \min_\bfw \big\{\Xi(\bfw) ~\big|~ \bfz_{\oT_*} = 0,~ \bfz_{T_*} \geq 0 \big. \big\},  \\
	& \eta_1 := \min_{\S \subsetneq T_*} \big\{  \min_\bfw \big\{\Xi(\bfw) ~\big|~ \bfz_{\omS} = 0,~ \bfz_\S \geq 0 \big. \big\} \big\}, \\	
	& \eta_2 := \min_{\S \nsubseteq T_*} \big\{  \min_\bfw \big\{\Xi(\bfw) ~\big|~ \bfz_{\omS} = 0,~ \bfz_\S \geq 0 \big. \big\} \big\}, 
\end{aligned}
\end{align}
where $\S \subsetneq T_*$ means $\S$ is a proper subset of $T_*$ and $\S \nsubseteq T_*$ means $\S$ is not a subset of $T_*$. Since Assumption \ref{assumption-slater+} holds, all the three values above are finite. It follows from Lem. \ref{lem-piecewise-global-phi+} that
\begin{equation} \label{piecewise-global-phi+-eta}
\min_{\bfw} \Xi(\bfw) + \Psi_{+,\bfmu}(\bfz) = \min \{ \eta_0, \eta_1, \eta_2 \}.
\end{equation}
In particular, $\eta_0$ is the optimal value of \eqref{convex-S} with $\S = T_*$. If we take $\S \subseteq T_*$, then the optimal value of \eqref{convex-S} must be greater than or equal to $\eta_0$. Noticing that we actually take $\S \subsetneq T_* $ in the definition of $\eta_1$ and the solution $\bfw^*$ satisfies \eqref{card-min}, then it holds that $\eta_1 - \eta_0 > 0$. 
%Moreover, it follows from Lem. \ref{lem-piecewise-global-phi+} that
%	\begin{equation} \label{piecewise-global-phi+-eta}
%	\min_{\bfw} \Xi(\bfw) + \Psi_{+,\bfmu}(\bfz) = \min \{ \eta_1, \eta_2, \eta_3 \}.
%\end{equation}
%Let us give some explanations about the above parameters. Firstly, $\eta_0$ is the optimal value of \eqref{omega+}. Secondly, if $\S \subseteq T_*$ holds, then the optimal value of the following convex program must be less or equal to $\eta_1$: 
%we consider the convex programs with the support sets being the subset of $\S$ and it holds that $\S \subsetneq T_*$
%In particular, \eqref{card-min} ensures that $\eta_1 - \eta_0 > 0$. The constant $\eta_2$ is finite when Assumption \ref{assumption-slater} holds.
\begin{theorem} \label{thm-D-local-global}
	Let us consider Problem \eqref{D} with $\Psi_\bfmu = \Psi_{+,\bfmu}$. Suppose that Assumption \ref{assumption-slater+} holds and $\bfw^*$ is a local minimizer satisfying \eqref{card-min}. If regularization parameter $\bfmu$ is taken as
	\begin{align} \label{mu}
		\left\{ \begin{aligned}
			& \sum_{i \in T_*} \mu_i \leq \eta_1 - \eta_0, \\
			& \mu_j > \max\{ \eta_1 - \eta_2, 0 \},~\mbox{for}~ j \in \oT_*,
		\end{aligned} \right.
	\end{align} 
	then $\bfw^*$ is a global minimizer of \eqref{D} with $\Psi_\bfmu = \Psi_{+,\bfmu}$. 
\end{theorem}
\bp
Since Lem \ref{lem-relation-D-D+} indicates that $\bfw^*$ is a global minimizer of \eqref{D+}, the first line in \eqref{mu} implies
\begin{align*}
	\Xi(\bfw^*) + \Psi_{+,\bfmu}(\bfz^*) = & \eta_0 + \sum_{i \in T_*} \mu_i \leq \eta_1 \notag \\
	\leq & \min_{\S \subsetneq T_*} \Big\{  \min_\bfw \Big\{\Xi(\bfw) + \sum_{i \in \S} \mu_i ~\Big|~ \bfz_{\omS} = 0, \bfz_\S \geq 0 \Big. \Big\} \Big\}.
\end{align*}
Furthermore, we can use \eqref{mu} to derive
\begin{align*}
	\Xi(\bfw^*) + \Psi_{+,\bfmu}(\bfz^*) = &\eta_0 + \sum_{i \in T_*} \mu_i \leq \eta_1 \leq \eta_2 + \min_{j \in \oT_*} \mu_j \\
	\leq & \min_{\S \nsubseteq T_*} \Big\{  \min_\bfw \Big\{\Xi(\bfw) + \sum_{i \in \S} \mu_i ~\Big|~ \bfz_{\omS} = 0,~ \bfz_\S \geq 0 \Big. \Big\} \Big\}, 
\end{align*}
where the first and second inequalities are from the first and second line of \eqref{mu} respectively. The last inequality above holds due to $\S \cap \oT \neq \emptyset$ when $\S \nsubseteq T^*$. Finally, using Lem. \ref{lem-piecewise-global-phi+}, we can conclude that $\bfw^*$ is a global minimizer of \eqref{D}.  
\ep

Next, let us consider Problem \eqref{D} with $\Psi_\bfmu = \Psi_{0,\bfmu}$. Lem. \ref{lem-relation-D-D0} indicates that $\bfw^*$ belongs to the following set of solution:
\begin{align*} 
	\Omega_0:= \argmin_{\bfw = [\bfy; \bfz] }\; 
	\Xi(\bfw) \quad \mbox{s.t.} \ \bfz_{\oT_*} = 0.
\end{align*}
We assume that $\bfw^*$ is the point in $\Omega_0$ with the smallest cardinality on $\bfz$ variable:  
\begin{align} \label{card-min0}
	\bfw^* \in \argmin_\bfw \{ \| \bfz \|_0 ~|~ \bfw \in \Omega_0 \}.
\end{align}
%If this does not hold, we can find the point in $\Omega^*$ with the smallest number of nonzero entries, and then reconstruct $\Omega^*$.
Similar to parameters $\eta_1$, $\eta_2$, and $\eta_3$, we define
\begin{align} \label{para-xi} 
	\begin{aligned} 
		& \xi_0 :=   \min_\bfw \big\{\Xi(\bfw) ~\big|~ \bfz_{\oT_*} = 0 \big. \big\},  \\
		& \xi_1 := \min_{\S \subsetneq T_*} \big\{  \min_\bfw \big\{\Xi(\bfw) ~\big|~ \bfz_{\omS} = 0 \big. \big\} \big\}, \\	
		& \xi_2 := \min_{\S \nsubseteq T_*} \big\{  \min_\bfw \big\{\Xi(\bfw) ~\big|~ \bfz_{\omS} = 0 \big. \big\} \big\}. 
	\end{aligned}
\end{align}
Since Assumption \ref{assumption-slater} holds, the above three values must be finite and it follows from Lem. \ref{lem-piecewise-global-phi0} that
\begin{equation} \label{piecewise-global-phi+-xi}
	\min_{\bfw} \Xi(\bfw) + \Psi_{0,\bfmu}(\bfz) = \min \{ \xi_0, \xi_1, \xi_2 \}.
\end{equation}
The constant $\xi_2$ is finite when Assumption \ref{assumption-slater} holds. Similar to the explanation on \eqref{para-eta}, we can derive $\xi_1 - \xi_0 > 0$. 
\begin{theorem} \label{thm-D-local-global0} Let us consider Problem \eqref{D} with $\Psi_\bfmu = \Psi_{0,\bfmu}$. Suppose that Assumption \ref{assumption-slater} holds and $\bfw^*$ is a local minimizer satisfying \eqref{card-min0}. If regularization parameter $\bfmu$ is taken as
	\begin{align} \label{mu0}
		\left\{ \begin{aligned}
			& \sum_{i \in T_*} \mu_i \leq \xi_1 - \xi_0, \\
			& \mu_j > \max\{ \xi_1 - \xi_2, 0 \},~\mbox{for}~ j \in \oT_*.
		\end{aligned} \right.
	\end{align} 
	then $\bfw^*$ is a global minimizer of \eqref{D} with $\Psi_\bfmu = \Psi_{0,\bfmu}$. 
\end{theorem}
The proof is similar to that of Thm. \ref{thm-D-local-global}. Therefore, we omit the details for brevity.
%\bp
%Since Lem \ref{lem-relation-D-D0} indicates that $\bfw^*$ is a global minimizer of \eqref{D0}, the first line in \eqref{mu} implies
%\begin{align*}
%	\Xi(\bfw^*) + \Psi_{0,\bfmu}(\bfz^*) = \eta_0 + \sum_{i \in T^*} \mu_i \leq \eta_1 \leq \min_{\S \subsetneq T^*} \Big\{  \min \Big\{\Xi(\bfw) + \sum_{i \in \S} \mu_i ~\Big|~ \bfz_{\omS} = 0 \Big. \Big\} \Big\}.
%\end{align*}
%Furthermore, we can use \eqref{mu} to derive
%\begin{align*}
%	\Xi(\bfw^*) + \Psi_{0,\bfmu}(\bfz^*) = \eta_0 + \sum_{i \in T^*} \mu_i \leq \eta_1 \leq \eta_2 + \min_{j \in \oT^*} \mu_j \leq \min_{\S \nsubseteq T^*} \Big\{  \min \Big\{\Xi(\bfw) + \sum_{i \in \S} \mu_i ~\Big|~ \bfz_{\omS} = 0 \Big. \Big\} \Big\}, 
%\end{align*}
%where the first and second inequalities are from the first and second line of \eqref{mu} respectively. The last inequality above holds due to $\S \cap \oT \neq \emptyset$ when $\S \nsubseteq T^*$. Finally, using Lem. \ref{lem-piecewise-global-phi0}, we can conclude that $\bfw^*$ must be a global minimizer of \eqref{D}.  
%\ep

\begin{remark} \label{remark-correspondence-global}
	Thm. \ref{thm-correspondence} establishes the correspondence between local minimizers of \eqref{P} and \eqref{D} for any given regularization parameters $\bfla$ and $\bfmu$. A more important question is the correspondence between the global solutions of these two problems. Thms. \ref{thm-D-local-global} and \ref{thm-D-local-global0} actually implies that this correspondence can be established provided that $\bfmu$ is appropriately selected. Generally, selecting this parameter is challenging because the index set $T_*$ is unknown in practice and parameters $\eta_i$ (resp. $\xi_i$) for $i = 0,1,2$ are difficult to compute. Nevertheless, the two theorems indicates that the weights on support set $T_*$ should be sufficiently small, whereas those on $\oT_*$ should be large enough. This observation may provide guidelines on numerical study. 
\end{remark}
\section{Conclusion} \label{Section-Conclusion}
This paper studies the stationary duality theory of CCOP. It extends the two-block model in \cite{zhang2025composite} to the three-block case with cardinality function $\Phi_\bfla \in \{ \Phi_{0,\bfla}, \Phi_{+,\bfla} \}$. The sufficient conditions for existence of global solutions of primal and dual problems are investigated. Compared with the existing literature of cardinality optimization, these conditions are easy to check and applicable to more general models. The one-to-one correspondence of solutions to the primal and dual problems is further established. These results lay theoretical foundation for solving the dual CCOP. It is noteworthy that the composite cardinality term has been reduced to a simple cardinality term in the dual CCOP. Therefore, the dual problem has a more favorable structure for algorithmic design and convergence analysis.
Such advantage has been demonstrated by the two-block case with $f$ being strongly convex in the our
previous work \cite{zhang2025composite}.
We will explore the three-block case with $f$ being just convex in the future work. 
%This requires more detailed investigation, and therefore we leave it for future work.

%
%\section*{Acknowledgment}
%This paper was supported by Hong Kong RGC General Research Fund (PolyU/15309223), PolyU AMA Projects (P0044200, P0045347), the National Natural Science Foundation of China (12571316, 12131004), the National Key R\&D Program of China (2023YFA1011100), and 111 Project of China (B16002).

%\section*{Acknowledgment}
%This paper was supported by Hong Kong RGC General Research Fund (PolyU/15309223), PolyU AMA Projects (P0044200, P0045347), the National Natural Science Foundation of China (12571316, 12131004), the National Key R\&D Program of China (2023YFA1011100), and 111 Project of China (B16002).

%\section*{Appendix}
%\begin{proposition}
%	If Assumption \ref{asm-coercive-f} holds, then $\liminf_{\| \bfx \| \to \infty} \Psi_\bfmu_k (\bfx) = \infty$.
%\end{proposition}

%%%%%%%%%%%%%%%%%%%%%%%%%%%%%%%%%%%%%%%%%%%%%%%%%%%%%%%%%%%%%%
\bibliographystyle{siamplain}
\bibliography{CCOPT_ref} % if more than one, comma separated

\end{document}